\documentclass[11pt]{article} 

\usepackage[utf8]{inputenc} 
\usepackage{authblk}
\usepackage{geometry}
\usepackage{float}		
\usepackage{comment}
\usepackage[english]{babel}
\usepackage{graphicx}
\usepackage{caption}
\usepackage{subcaption}
\usepackage{verbatim}
\usepackage{amsmath}
\usepackage{amssymb}
\usepackage{amsfonts}
\usepackage{ulem}
\usepackage{bbold}

\usepackage{amsthm}

\theoremstyle{plain}

\newtheorem{thm}{Theorem}

\newtheorem{lem}[thm]{Lemma}
\newtheorem*{lem*}{Lemma}
\newtheorem{prop}[thm]{Proposition}
\newtheorem*{prop*}{Proposition}

\newtheorem{cor}[thm]{Corollary}
\theoremstyle{definition}

\theoremstyle{remark}

\newtheorem{rem}{Remark}

\theoremstyle{plain}
\newtheorem{hyp}{Hypothesis}

\numberwithin{equation}{section}

\newcommand\RR{\mathbb{R}} 

\usepackage{cancel}

\def\restriction#1#2{\mathchoice
	{\setbox1\hbox{${\displaystyle #1}_{\scriptstyle #2}$}
		\restrictionaux{#1}{#2}}
	{\setbox1\hbox{${\textstyle #1}_{\scriptstyle #2}$}
		\restrictionaux{#1}{#2}}
	{\setbox1\hbox{${\scriptstyle #1}_{\scriptscriptstyle #2}$}
		\restrictionaux{#1}{#2}}
	{\setbox1\hbox{${\scriptscriptstyle #1}_{\scriptscriptstyle #2}$}
		\restrictionaux{#1}{#2}}}
\def\restrictionaux#1#2{{#1\,\smash{\vrule height .8\ht1 depth .85\dp1}}_{\,#2}}

\usepackage{color}

\newcommand\p\varphi
\newcommand\e\varepsilon

\newcommand{\cred}[1]{{#1}}
\newcommand{\am}[1]{{#1}}

\graphicspath{{./figures_CGM/}}

\title{Limiting Hamilton-Jacobi Equation for the Large Scale Asymptotics of a Subdiffusion Jump-Renewal Equation}

\author[1]{Calvez, Vincent \thanks{vincent.calvez@mat.cnrs.fr}}
\author[2]{Gabriel, Pierre\thanks{pierre.gabriel@uvsq.fr}}
\author[3,4,5]{Mateos González, Álvaro \thanks{alvaro.mateos-gonzalez@umontpellier.fr}}
\affil[1]{Institut Camille Jordan, {Université Claude Bernard Lyon 1},
	43 boulevard du 11 novembre 1918, F-69622 Villeurbanne Cedex, {France}}
\affil[2]{Laboratoire de Mathématiques de Versailles, {UVSQ, CNRS, Université Paris-Saclay},
	45 avenue des États-Unis,
	F-78035 Versailles cedex, {France}}
\affil[3]{Institut Montpellli\'erain Alexander Grothendieck (IMAG), UMR CNRS 5149, {Universit\'e de Montpellier},
	34090 Montpellier, {France}}
\affil[4]{Institut des Sciences de l'\'Evolution de Montpellier (ISEM), UMR CNRS 5554, {Universit\'e de Montpellier},
	34095 Montpellier, {France}}
\affil[5]{MISTEA, UMR CNRS 0729, {INRA and SupAgro Montpellier},
	34060 Montpellier, {France}}

\begin{document}

\maketitle

\begin{abstract}
Subdiffusive motion takes place at a much slower timescale than diffusive motion. As a preliminary step to studying reaction-subdiffusion pulled fronts, we consider here the hyperbolic limit $(t,x) \to (t/\e, x/\e)$ of an age-structured equation describing the subdiffusive motion of, \emph{e.g.}, some protein inside a biological cell. Solutions of the rescaled equations are known to satisfy a Hamilton-Jacobi equation in the formal limit $\e \to 0$. In this work we derive uniform Lipschitz estimates, and establish the convergence towards the viscosity solution of the limiting Hamilton-Jacobi equation. The two main obstacles overcome in this work are the non-existence of an integrable stationary measure, and the importance of memory terms in subdiffusion.	
\end{abstract}

\vspace{0.5cm}

\small{\noindent
}

\vspace{0.5cm}




	\section{Introduction}
	\label{CGM_sec_intro}

	\subsection{Model description} \label{CGM_ss_model}
	
	Consistent experimental evidence stemming from recent methodological advances in cell biology such as \textit{in vivo} single molecule tracking, report that the intra-cellular random motion of certain molecules often deviates from Brownian motion. Macroscopically, their mean squared displacement does not scale linearly with time, but as a power law $t^\mu$ for some exponent $0 < \mu < 1$~\cite{Cox,Bronstein,Parry2014,DiRienzo2014,Izeddin}. This behaviour, due to crowding and trapping phenomena, is usually referred to as `anomalous' diffusion or `subdiffusion'. The reader may consult \cite{Hoefling-Franosch} for a review.
	
	One of the standard mechanisms used to describe the emergence of subdiffusion in cells is continuous time random walks (CTRW), a generalisation of random walks that couples a waiting time random process at each `jump' of the random walk \cite{Montroll1965}. CTRW can be used \cite{Metzler2000, Metzler2004, Reaction-Transport} to derive macroscopic equations governing the spatiotemporal dynamics of the density of random walkers located at position $x$ at time $t$:
	\[
	\partial_t\rho(x,t)=D_\mu \mathcal{D}_t^{1-\mu}\Delta\rho(x,t) .
	\]
	Here, $D_\mu$ is a generalised diffusion coefficient and
	$\mathcal{D}_t^{1-\mu} (f)(t) = \frac 1 {\Gamma(\mu)} \frac{{\rm d}}{{\rm d}t} \int_0^t \frac{f(t')}{(t-t')^{1-\mu}} {\rm d}t'$ is the Riemann-Liouville fractional derivative operator. Such a fractional dynamics formulation is very attractive for modelling in biology, in particular because of its apparent similarity with the classical diffusion equation. However, contrary to the diffusion equation, the Riemann-Liouville operator is non-local in time. This is the `trace' of the non-Markovian property of the underlying CTRW process. Indeed, memory terms play a crucial role in subdiffusive processes. This non-Markovian property becomes a serious obstacle when one wants to couple subdiffusion with chemical reaction \cite{Henry2006,Yuste2008,Fedotov2014}.  
	
	In this work, following~\cite{VR_2002},
	we take an alternative approach that rescues the Markovian property of the jump process at the expense of a supplementary age variable. We associate each random walker with a residence time (age, in short) $a$, which is reset when the random walker jumps to another location. We denote by $n(t,x,a)$ the probability density function of walkers at time $t$ that have been located at $x$ exactly during the last span of time $a$. The dynamics of the CTRW are then described~\cite{VR_2002, YH_2006, Reaction-Transport, Falconer2012} by means of an age-renewal equation with spatial jumps:
	\begin{equation}
	\left\{
	\begin{aligned}
	\label{CGM_eq_n}
	& \partial_t n(t,x,a) + \partial_a n(t,x,a) + \beta(a) n(t,x,a) = 0\, ,\quad t \geq0, \quad a > 0\, , \quad x\in \RR^d\smallskip\\
	& n(t,x,a=0) = \int_0^\infty \beta(a') \int_{\RR^d} \omega(x-x') n(t,x',a') \,{\rm d}x' \,{\rm d}a' \smallskip\\
	& n(t = 0,x,a) = n^0(x,a). 
	\end{aligned}\right.
	\end{equation}
	The boundary condition on $n(t,x,0)$ at age $a = 0$ accounts for the particles landing at position~$x$ at time~$t$ after having `jumped' from position $x'$, at which they had remained during a time span exactly equal to $a'$. Here, $\beta$ is the age-dependent  rate of jump, and $\omega$ is the distribution of jump distances. They are chosen in the following way.

	\begin{hyp}[Space jump kernel $\omega$ and jump rate $\beta$] \label{CGM_hyp_omega}~\\
		We assume that $\omega$ is an isotropic multivariate normal distribution of mean $0$ and variance $\sigma^2$, and  that $\beta$ is decaying for large age in a precise way: 
		\begin{equation} \label{CGM_eq_omega_beta}
		\left\{
		\begin{aligned}
		& \omega(x) = \frac 1 {(\sigma \sqrt{2 \pi})^d} \exp \left( - \frac {|x|^2}{2 \sigma^2}\right) , & \sigma > 0  \\
		& \beta(a)  =  \frac \mu {1+a},    & 0 < \mu < 1.
		\end{aligned} \right.
		\end{equation}
		The assumption of a Gaussian $\omega$ \am{can be relaxed to even functions that exhibit an exponential decay faster than the decay of the initial condition, as stated in Corollary~\ref{cor_exp_tails}.} 
		However, the normal distribution provides simpler asymptotic estimates on the Hamiltonian $H$ in Section~\ref{sec_prop_H} and allows proofs to be clearer.

		
		 
		 The specific choice of the rate of jump $\beta$ is crucial. Only the case $\mu \in (0,1)$ yields subdiffusion. This could be relaxed to an asymptotic equivalence -- as is the case in the study of the related space-homogeneous problem in~\cite{BLM} --, however we will stick to \eqref{CGM_eq_omega_beta} for the sake of clarity.
	\end{hyp}

	The fact that the loss term $\beta(a) n(t,x,a)$ is recovered in the boundary condition (and that $\omega$ is a probability distribution) leads to the conservation of the total population density $\int_0^\infty \int_{\RR^d} n(\cdot, x, a) \,{\rm d}x \,{\rm d}a$ along time.

	We restrict to initial conditions compactly supported in age. \cred{More precisely we have the following assumption:}
	\begin{equation}\label{eq:supp compact}
	(\forall x)\quad {\rm supp} (n^0(x,\cdot)) = [0,1] .
	\end{equation}
	Further technical hypotheses will be made later on.

	The probability that a particle reaches age $a$ without jumping is $\exp \left(- \int_{0}^a \beta(s) \,{\rm d}s \right)$. On the other hand, the jump rate of particles at age $a$ is $\beta(a)$. Hence, the distribution of residence times $\Phi(a)$ (meaning the distribution of the age of particles when they jump) is given by
	\begin{equation} \label{CGM_eq_def_captial_Phi}
	\Phi(a) = \beta(a) \exp \left( -\int_0^a  \beta(s) \,{\rm d}s \right) = \frac \mu {(1+a)^{1+\mu}}.
	\end{equation}
	A noteworthy observation is that the mean residence time of particles $\int_0^\infty a \Phi(a)\, da$ is infinite since $\mu \in (0,1)$. This is a signature of subdiffusion at a larger scale \cite{Reaction-Transport}.

\cred{Our motivation is the asymptotic analysis of pulled fronts in reaction-subdiffusion equations in the hyperbolic regime $(t / \e, x/ \e, a)$. On the one hand, reaction-subdiffusion equations have stimulated an extensive literature \cite{fronts_ABAA_2008, fronts_ABAA_2010, VNN_2013, fronts_Nepo_Volp_2013, fronts_Polish_2013, N_2015, fronts_time_scales}. On the other hand, classical pulled reaction-diffusion fronts have been studied in the same hyperbolic regime $(t/\e,x/\e)$ by means of stochastic calculus methods [Freidlin SIAM J APPL MATH 1986] and PDE methods [Evans-Souganidis Indiana J. 1989]. The singular limit yields a Hamilton-Jacobi equation that encodes the motion of the level set of the solution. Here, we extend rigorously this analysis for the subdiffusion equation \eqref{CGM_eq_n} in the absence of reaction.} 


	\subsection{Hyperbolic limit and derivation of the Hamilton-Jacobi equation.}
	
	\cred{We perform the Hopf-Cole transformation in order to study the large scale asymptotics:}
	\begin{equation} \label{CGM_eq_def_phi}
	n_\e (t,x,a) = n \left( t / \e, x/\e, a \right) = \exp \left( - \phi_\e(t,x,a) / \e \right).
	\end{equation}
	This enables us to accurately measure the behaviour of small, exponential tails of the probability density function $n$, reminiscent of large deviation principle theory. 

The function $n_\e$ satisfies the following equation,
	\begin{equation} \label{CGM_eq_n_e}
	\left\{
	\begin{aligned}
	& \partial_t n_\e + \frac 1 \e \partial_a n_\e + \frac 1 \e \beta n_\e = 0\, ,\quad t \geq0, \quad a > 0\, , \quad x\in {\RR^d}\smallskip\\
	& n_\e(t,x,0) = \int_0^{1+ t / \e} \int_{\RR^d} \beta(a) \omega(z) n_\e(t,x - \e z ,a) \,{\rm d}z \,{\rm d}a \smallskip\\
	&	n_\e(0,x,a) = n_\e^0(x,a) = n^0(x/\e, a), 
	\end{aligned}\right.
	\end{equation}
	\am{where the upper integration bound $1+t/\e$ is the upper bound of the support in age of $n_\e$ at time $t$, due to the transport of the compact support of the initial condition.}
	
\cred{Accordingly, the function} $\phi_\e$ satisfies the following \cred{non-linear problem}:
	\begin{equation}
	\label{CGM_eq_phi}
	\left\{
	\begin{aligned} 
	&	\partial_t \phi_\e + \frac 1 \e \partial_a \phi_\e - \beta = 0\, ,\quad t \geq0, \quad a > 0\, , \quad x\in \RR^d \smallskip\\
	&	\exp\left(-\phi_\e(t,x,0) / \e \right) = \int_0^{1+t/\e} \int_{\RR^d} \beta(a) \omega(z) \exp \left( - \phi_\e(t,x - \e z ,a) / \e\right) \,{\rm d}z \,{\rm d}a \smallskip\\
	&	\phi_\e(0,x,a) = \phi_\e^0(x,a) = -\e \ln \left(n^0(x/\e, a) \right).
	\end{aligned} \right.
	\end{equation}
	Let us denote by $\psi_\e$ the boundary value at $a=0$, which will be our main unknown:
	\begin{equation}
	\psi_\e(t,x) = \phi_\e(t,x,0).
	\end{equation}
	We compute the solution of equation \eqref{CGM_eq_phi} along characteristic lines:
	\begin{equation} \label{CGM_eq_carac}
	\phi_\e (t,x,a) =  \begin{cases}
	\psi_\e(t-\e a, x) + \e \int_0^a \beta (s) \,{\rm d} s, & \qquad  t > 0, \, \varepsilon a < t \smallskip\\
	\phi_\e^0( x, a - t / \e) + \e \int_{a-t / \e}^a \beta (s) \,{\rm d} s, & \qquad  t \geq 0, \, a \geq t / \e .
	\end{cases}
	\end{equation}
	We inject~\eqref{CGM_eq_carac} into the second line of~\eqref{CGM_eq_phi} \cred{so as to get} 
	\begin{multline} \label{CGM_eq_eq}
	1 = \int_0^{t/\e} \Phi(a) \int_{\RR^d} \omega(z)  \exp\left( \frac 1 \e \left[ \psi_\e (t,x) - \psi_\e(t-\e a, x - \e z) \right] \right) \,{\rm d} z \,{\rm d} a \\
	+ \int_{t/\e}^{1+t/\e} \Phi(a) \int_{\RR^d} \omega(z)  \exp\left(\frac 1 \e \left[ \psi_\e (t,x) - \phi_\e^0(x - \e z, a - t / \e) \right] + \int_0^{a - t/\e} \beta \right) \,{\rm d} z \,{\rm d} a.
	\end{multline}
	Taking the formal limit of \eqref{CGM_eq_eq} when $\varepsilon \to 0$ yields \cred{the following Hamilton-Jacobi equation}:
	\begin{equation}
	\label{CGM_eq_limiting_HJ}
	1 = \int_0^\infty \Phi(a) \exp\left(a \partial_t \psi_0 (t,x)\right) \,{\rm d}a \int_{\RR^d} \omega(z) \exp\left(z \cdot \nabla_x \psi_0 (t,x)\right) \,{\rm d} z .
	\end{equation}
	\cred{Observe that it} is equivalent to:
	\begin{equation}
	\label{CGM_HJ}
	\partial_t \psi_0 (t,x) + H(\nabla_x \psi_0)(t,x) = 0 ,
	\end{equation}
	with $H$ defined as follows, where $\hat{\Phi}^{-1}$ is the inverse function of the Laplace transform of $\Phi$: 
	\begin{equation} \label{CGM_eq_def_H}
	H(p) = \hat{\Phi}^{-1}\left( \frac 1 { \int_{\RR^d} \omega(z) \exp(z \cdot p) \,{\rm d}z } \right) .
	\end{equation}
	
	\begin{rem}[\cred{About the scaling}] 
\cred{We emphasize that the limiting equation \eqref{CGM_HJ} makes sense for a large class of functions $\beta$, including constant  rates of jump. 
	On the contrary, diffusion limits depend on the decay properties of $\beta$, as illustrated by the anomalous scaling   $\left( \frac t {\e^{2/\mu}}, \frac x \e, a \right)$ under which they are performed~\cite{Reaction-Transport}. \am{In our scaling, the slow decay of $\beta$ has an impact on the properties of the Hamiltonian function $H$ and also on the estimates that we are able to derive in the proof of convergence}.  
} 
	\end{rem}

	We discuss several properties of the Hamiltonian $H$ in Section~\ref{sec_prop_H}: its smoothness, coercivity, convexity but lack of strict uniform convexity, and its asymptotic behaviour near $0$ and $\infty$.
	
	We recall that, under suitable hypotheses on the Hamiltonian $H$ and on the initial condition $g$, classical existence and uniqueness results hold for the evolution Hamilton-Jacobi Cauchy problem:
	\begin{equation}\label{CGM_HJ_CP}
	\left\{ \begin{aligned}
	\partial_t u(t,x) + H(\nabla_x u(t,x)) = 0, & \qquad (t,x) \in (0,T) \times \RR^d \\
	u(0,x) = g(x), & \qquad x \in \RR^d.
	\end{aligned} \right.
	\end{equation}
	
	We state hereafter a relevant uniqueness theorem in a suitable class of functions: a version of~\cite[Theorems 19.11 and 19.17]{book_Clarke} for a homogeneous Hamiltonian that is not polynomially bounded above.
	
	\begin{thm}[Uniqueness theorem]\label{CGM_thm_uniqueness}
		Let $H$ be locally Lipschitz, convex and superlinear. Let $g$ be bounded below and Lipschitz continuous.
		Then there exists a unique viscosity solution of~\eqref{CGM_HJ_CP} within the class of Lipschitz continuous functions.
	\end{thm}

	This uniqueness theorem is a corollary of~\cite[Corollary~19.17]{book_Clarke}, which follows from~\cite[Theorem~19.11]{book_Clarke}. In that last theorem it is assumed that $H$ has polynomial growth for $|p| \to \infty$, which is not our case, as stated in Proposition~\ref{CGM_prop_H_at_infty}. We overcome this issue by assuming that $u$ is globally Lipschitz continuous so that $H$ \cred{can be restricted to a compact set.}  

	\subsection{Main hypotheses and results} \label{CGM_ssec_plan}
	
	

In this work we establish the rigorous proof of convergence from \eqref{CGM_eq_phi} to \eqref{CGM_HJ} as $\e \to 0$, under suitable hypotheses on the initial data. 
	

	
	\begin{hyp}[Initial condition $\phi_\e^0$] \label{CGM_hyp_IC}~\\
		We assume that the  initial condition has the following form:
		\begin{equation}\label{CGM_eq_Ansatz}
		\phi_\e^0 (x,a) = v_\e(x) + \e \eta_\e(x,a) + \chi_{[0,1]}(a)\,,
		\end{equation}
		where $\eta_\e(x,0) = 0$ by convention.
		Here, $\chi$ denotes the convex characteristic function: $\chi_{[0,1]} (a) = 0$ for $a \in [0,1]$ and $+\infty$ for $a > 1$. Hence, $\phi_\e^0$ takes finite value in $[0,1]$ only, according to the assumption on $n^0$ \eqref{eq:supp compact}. The functions $\phi_\e^0$, $v_\e$ and $\eta_\e$ satisfy the following properties uniformly over $\varepsilon$:
		\begin{enumerate}
			\item \label{CGM_H2_sublinearity} \label{CGM_hyp_IC_v} $v_\e$ is bounded below. 
			
			\item \label{CGM_H2_eta} \label{CGM_hyp_IC_eta} \label{CGM_hyp_IC_inf_eta}  \label{CGM_hyp_IC_sup_eta}  $\eta_\e$ is bounded uniformly in $\e$.
			
			\item \label{CGM_H2_Lips} \label{CGM_hyp_IC_regularity} $\phi_\e^0$ is Lipschitz continuous in $x$ uniformly in $a \in [0,1]$: there exists $C_L$ such that, for any $a \in [0,1]$, for any $x_1, x_2 \in \RR^d$ and for any $\e > 0$,
			\begin{equation}\label{CGM_eq_Lips}
			\left| \phi_\e^0 (x_1,a) - \phi_\e^0 (x_2,a) \right| \leq C_L | x_1 - x_2 |.  
			\end{equation}
			
			\item \label{CGM_H2_semiconcavity} \label{CGM_hyp_IC_semi_concavity} $\phi_\e^0$ is semi-concave in $x$ uniformly in $a$: there exists $\mathfrak C_{xx} \in \RR$ such that for any $x, h \in \RR^d$ and $a \in [0,1]$, for any $\e > 0$,
			\begin{equation}
			\phi_\e^0 (x + h, a) + \phi_\e^0 (x - h, a) - 2 \phi_\e^0 (x, a) \leq \mathfrak C_{xx} |h|^2, \\
			\end{equation}
			(Or equivalently, $x \mapsto \frac{\mathfrak{C}_{xx}} {2} |x|^2 - \phi_\e^0(x, a)$ is convex, or  \allowbreak \cred{\mbox{$D_x^2 \phi_\e^0 \leq \mathfrak C_{xx}$} in the sense of distributions.})

			\item We assume that there exists a limit function $v$ such that $v_\e \xrightarrow[\e \to 0]{} v$, locally uniformly in $x, a$.
			
		\end{enumerate}
	\end{hyp}

	The following theorem is our main result.
	\begin{thm}\label{CGM_thm_viscosity_solution}
		Under Hypotheses \ref{CGM_hyp_omega} and \ref{CGM_hyp_IC}, $\psi_\e \xrightarrow[\e \to 0]{L^\infty_{\rm{loc}}} \psi_0$, which is the unique viscosity solution of the limiting Hamilton-Jacobi equation~\eqref{CGM_eq_limiting_HJ} with initial condition $v(x)$ among the class of bounded below, Lipschitz continuous functions. 
	\end{thm}

\cred{	The reader will find in Appendix~\ref{ssec_motivations} a comprehensive discussion about the hypotheses and some highlights of the proof.}

	\begin{rem}[Initial conditions -- interpretation] \label{rmk_IC_interpretation}
	\cred{The initial condition has the following shape in the original unknown:} 
		\[
		n_\e^0(x,a)=\tilde n_\e(x,a) \exp(-v_\e(x)/\e) \mathbb{1}_{[0,1]}(a),
		\]
\cred{So, for technical reasons, we restrict the age support to be uniformly bounded. Moreover, the initial profile is assumed to be  uniformly bounded below, locally in space, see \ref{CGM_ssec_persistence} for a discussion.}
%
	\end{rem}
	\subsection{Organization of the article}
	
	Section~\ref{CGM_sec_bounds} deals with the regularity of the solution which in turm yields compactness of $(\psi_\e)_\e$. 
	In Section~\ref{CGM_sec_viscosity_solution} it is established that $\psi_0$ is the unique viscosity solution of the limiting Hamilton-Jacobi equation~\eqref{CGM_eq_limiting_HJ}.
	
	During the first revision stage of this manuscript, the authors became aware of a preprint by Nordmann, Perthame, and Taing - now published~\cite{NPT2017}, which adresses similar questions in the context of evolutionary biology. Our model is simpler as it is conservative, and jump rates are homogeneous with respect to the space variable. On the other hand, our results are stronger as we \cred{establish} the rigorous limit of the problem as $\varepsilon\to 0$.

\section{Properties of the Hamiltonian}\label{sec_prop_H}

We will now prove that the Hamiltonian $H$ satisfies some properties often encountered in the literature.

\begin{prop} \label{CGM_prop_H}
	The Hamiltonian $H$ defined in \eqref{CGM_eq_def_H} has the following properties:
	\begin{enumerate}
		\item[(i)]
		$H \in C^\infty(\RR^d, \RR_+)$,
		\item[(ii)]
		$H$ exhibits quadratic growth at infinity.
		\item[(iii)]
		H is convex, but not strictly uniformly convex.
	\end{enumerate}

\end{prop}
\begin{proof}~
	\begin{enumerate}
		\item[(i)] Let
		\[
		F(p,h) = \int_0^\infty \Phi(a) \exp \left(- a h \right) \,{\rm d}a - \left( \int_{\RR^d} \omega(z) \exp \left(z \cdot p \right) \,{\rm d}z \right)^{-1}.
		\]
		$F$ is strictly decreasing with respect to its second variable over $\RR_+$. For all $p \in \RR^d \setminus \{0\}$, since $\omega$ is an isotropic multivariate Gaussian distribution centred at $0$ and $\Phi$ is a probability measure, it follows that $F(p,0) > 0$. For any $p \in \RR^d$, we have $\lim_\infty F(p, \cdot) < 0$. Hence for each $p \in \RR^d$ there exists a unique $H \in \RR_+$ such that $F(p,H) = 0$. This condition is equivalent to equation~\eqref{CGM_eq_limiting_HJ}, hence $H$ is well defined.
		
		The function $F$ is $C^\infty$, and $F(0,0) = 0$. Strict monotonicity and the implicit function Theorem yield the proof.
		\item[(ii)]
		We have:
		\begin{equation}\label{CGM_eq_aux_omega}
		\begin{aligned}
		\int_{\RR^d} \omega(z) e^{z \cdot p} \,{\rm d} z = \exp \left(\frac{\sigma^2 |p|^2}{2}\right).
		\end{aligned}
		\end{equation}
		It follows from equation~\eqref{CGM_eq_limiting_HJ} that
		\[
		\int_0^\infty \Phi(a) e^{-a H(p)} \,{\rm d}a \geq \int_0^1 \Phi(a) e^{- H(p)} \,{\rm d}a \geq C e^{-H(p)},
		\]
		hence
		\[
		C e^{-H(p)} \leq \exp\left(- \frac{ \sigma^2 |p|^2}{2}\right),
		\]
		which implies $H(p) \gtrsim 1+|p|^2$. Hence $H$ is coercive.
		\item[(iii)]
		Differentiating equation~\eqref{CGM_eq_limiting_HJ} with respect to $p$ yields the following \cred{identity}: 
		\[
		0 = \int_0^\infty \int_{\RR^d} \left( \nabla_p H(p) - \frac z a \right) a \Phi(a) \exp\left(-a H(p)\right) \omega(z) \exp\left(z \cdot p\right) \,{\rm d}z\,{\rm d}a.
		\]
\cred{Another step of differentiation gives us:}		
		\[
		0 = \int_0^\infty \int_{\RR^d} a D^2_p H(p) d\gamma(z,a) - \int_0^\infty \int_{\RR^d} a^2 \left( \nabla_p H (p) - \frac z a \right) \cdot \left( \nabla_p H (p) - \frac z a \right)^T d\gamma(z,a),
		\]
		where $d\gamma(z,a) = \Phi(a) \omega(z) \exp\left(- a H(p)\right) \exp\left(z \cdot p\right) \,{\rm d}z \,{\rm d}a$ is a non-negative measure. For any $z \in \RR^d$ and $a>0$, the  matrix $M := \left( \nabla_p H (p) - \frac z a \right) \cdot \left( \nabla_p H (p) - \frac z a \right)^T$ is symmetric, hence for all $h \in \RR^d$, $h^T M h \geq 0$. By integration, it follows that $D_p^2 H$ 
		is positive semi-definite.
		
		However, the Hamiltonian $H$ is not strictly uniformly convex, since $D^2_p H(0) = 0$. This is proved as follows. We first remark, from \eqref{CGM_eq_def_H} and \eqref{CGM_eq_limiting_HJ}, that $H(0) = 0$. We recover the following expression for $\nabla_p H(p)$:
		\[
		\nabla_p H(p) = \frac{ \int_0^\infty \int_\RR z \Phi(a) \omega(z) \exp(z \cdot p) \,{\rm d}z \,{\rm d}a}{\int_0^\infty \int_\RR a \Phi(a) \omega(z) \exp(z \cdot p) \,{\rm d}z \,{\rm d}a}.
		\] 
		Since $\int_0^\infty a \Phi(a) \,{\rm d}a = \infty$, we deduce that $\nabla_p H(0) = 0$ and recover:
		\[
		D^2_p H(0) = \frac {\int_\RR z \cdot z^T \omega(z) \,{\rm d} z}{\int_0^\infty a \Phi(a) \,{\rm d} a} = 0.
		\]
		
	\end{enumerate}
\end{proof}

\begin{prop}[Behaviour of $H$ around $0$] \label{CGM_prop_H_at_0}
	Around $|p| = 0$, we have 
	\begin{equation} \label{CGM_eq_H_at_0}
	H(p) \sim_{0} (\sigma |p|)^{2/\mu} \left( 2 \Gamma(1-\mu) \right)^{1/\mu}.
	\end{equation}
\end{prop}
\begin{proof}
	We have $\Phi(a) = \mu (1+a)^{-1-\mu}$, hence, thanks to equations~\eqref{CGM_eq_limiting_HJ} and~\eqref{CGM_eq_aux_omega}:
	\[
	\int_0^\infty \frac \mu {(1+a)^{1+\mu}} \left( \exp(-aH) - 1 \right) \,{\rm{d}}a = \hat \Phi (H) - 1 = \exp\left(-(\sigma |p|)^2 / 2 \right) - 1 \sim_{|p| = 0} - (\sigma |p|)^2 / 2.
	\]
	Denoting $b = aH$, since $H(0) = 0$ the left hand side becomes:
	\begin{align*}
	& \frac 1 H \int_0^\infty \frac \mu { \left( 1 + b/H \right)^{1+\mu}} (e^{-b} - 1) \,{\rm d}b \\
	= & H^\mu \int_0^\infty \frac \mu { \left( H/b + 1 \right)^{1+\mu} } b^{-1-\mu} (e^{-b} - 1) \,{\rm d}b \\
	\sim_{H=0} \; & H^\mu \int_0^\infty \mu b^{-1-\mu} \left( e^{-b} - 1 \right) \,{\rm d}b.
	\end{align*}
	Integrating that last expression by parts ends the proof.
\end{proof}

\begin{prop}[Behaviour of $H$ for large $|p|$]\label{CGM_prop_H_at_infty}
	Around $\infty$, we have 
	\begin{equation} \label{CGM_eq_H_at_infty}
	H(p) \sim_{\infty} \mu \exp \left( \frac{\sigma^2 |p|^2}{2}\right).
	\end{equation}
\end{prop}

\begin{proof}
	\cred{Back to} the computations of Proposition~\ref{CGM_prop_H_at_0}, we find:
	\begin{equation*}
	\frac{1}{H(p)} \int_0^\infty \frac{\mu}{ \left(1 + \frac{b}{H(p)}\right)^{1+\mu}} e^{-b} \,{\rm d}b = \exp \left( - \frac{\sigma^2 |p|^2}{2}\right).
	\end{equation*}
	The divergence of $H$ (\emph{e.g.} coercivity in Proposition~\ref{CGM_prop_H}) and the exponential tail $e^{-b}$ allow us to give the following equivalent: 
	\begin{equation*}
	\frac{1}{H(p)} \int_0^\infty \mu \left[ 1 - (1+\mu) \frac b {H(p)} \right]  e^{-b} \,{\rm d}b \sim_\infty \exp \left( - \frac{\sigma^2 |p|^2}{2}\right), 
	\end{equation*}
	which by integration leads to
	\begin{align*}
	\frac{\mu}{H(p)} \left(1 - \frac{1+\mu}{H(p)}\right) \sim_\infty 
	\exp \left( - \frac{\sigma^2 |p|^2}{2}\right) \\
	H(p) 
	\sim_\infty \mu \exp \left( \frac{\sigma^2 |p|^2}{2}\right) \left(1 - \frac{1+\mu}{H(p)}\right).
	\end{align*}
	The \cred{limit $H(p)\to +\infty$ as $|p|\to +\infty$} concludes the proof.
	
\end{proof}

\begin{figure}
	\begin{minipage}{0.47\textwidth}
		\centering{\includegraphics[width=\linewidth]{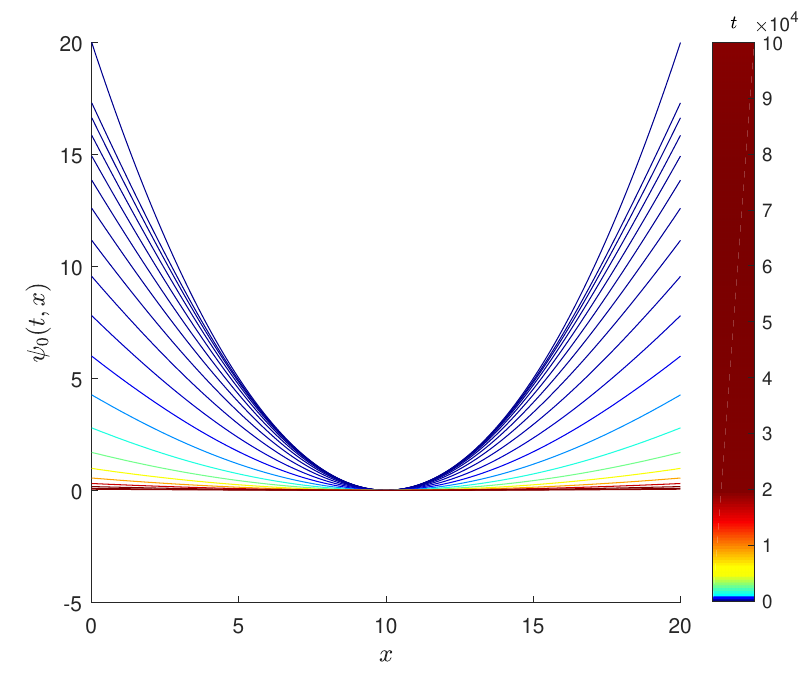}} 
	\end{minipage}
	\begin{minipage}{0.47\textwidth}
		\centering{\includegraphics[width=\linewidth]{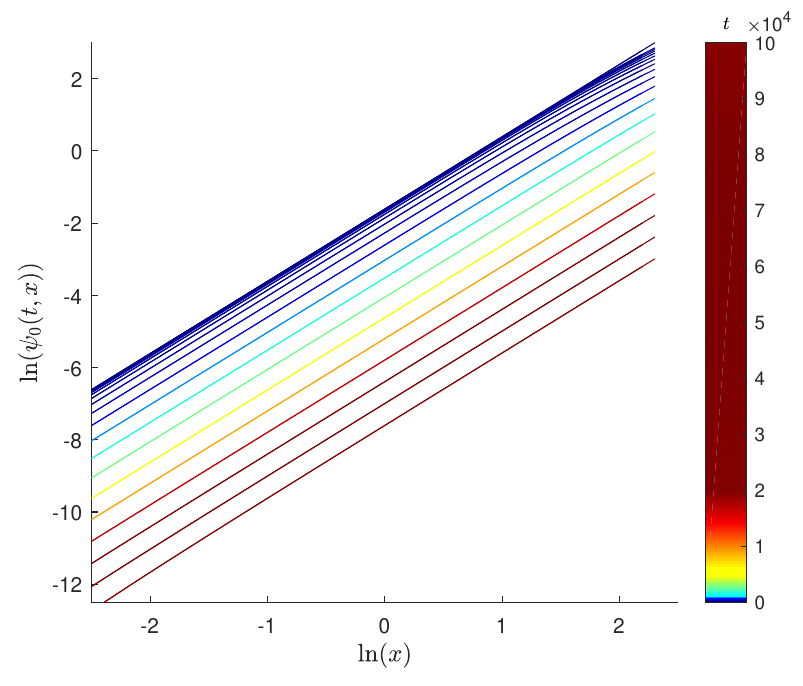}} 
	\end{minipage}
	{\begin{center}{\vspace{-1.1em} Diffusive case with $D = 0.01$ and $\beta (a) = D$. $\psi_0(0,x) = 0.2 (x-10)^2$. \vspace{-0.4em}}\end{center}}
	
	\begin{minipage}{0.47\textwidth}
		\centering{\includegraphics[width=\linewidth]{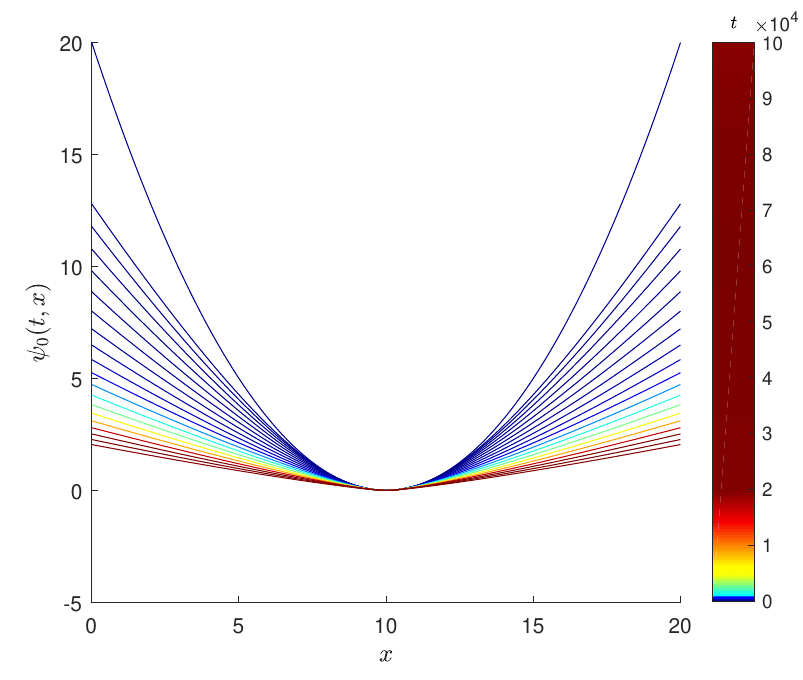}}
	\end{minipage}
	\begin{minipage}{0.47\textwidth}
		\centering{\includegraphics[width=\linewidth]{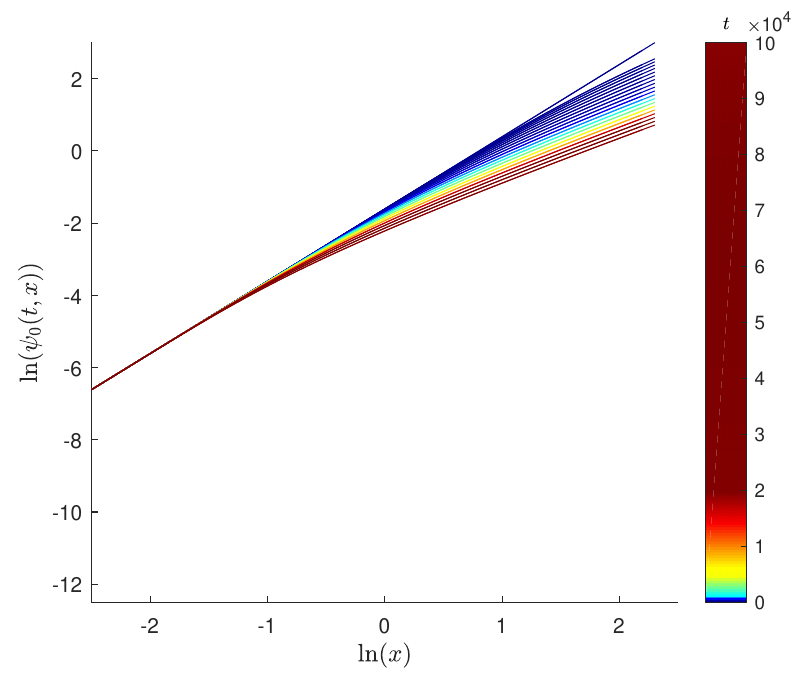}}
	\end{minipage}
	{\begin{center}{\vspace{-1.1em} Subdiffusive case with $\mu = 0.3$ and $\beta (a) = \mu / (1+a)$. $\psi_0(0,x) = 0.2 (x-10)^2$. \vspace{-0.2em}}\end{center}}
	
	\begin{minipage}{0.47\textwidth}
		\centering{\includegraphics[width=\linewidth]{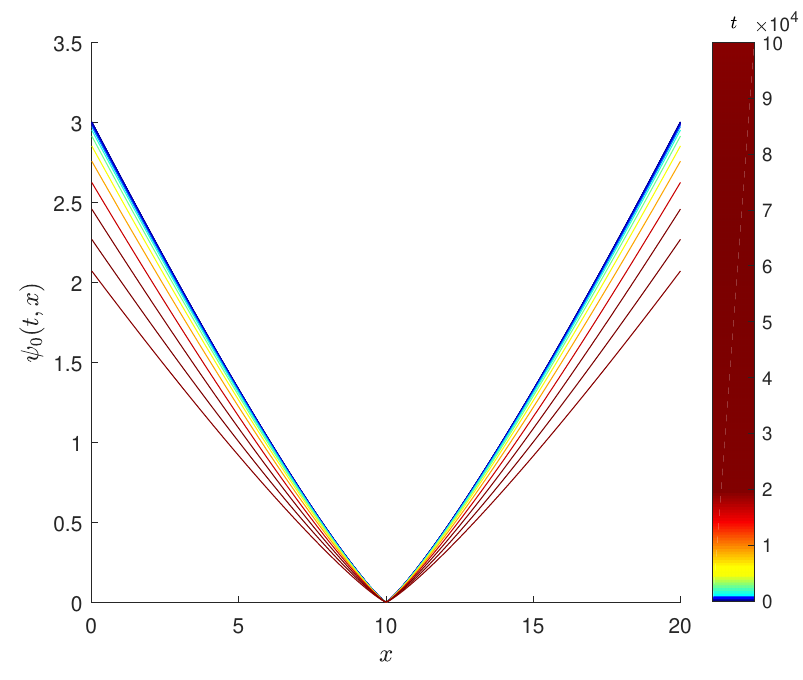}} 
	\end{minipage}
	\begin{minipage}{0.47\textwidth}
		\centering{\includegraphics[width=\linewidth]{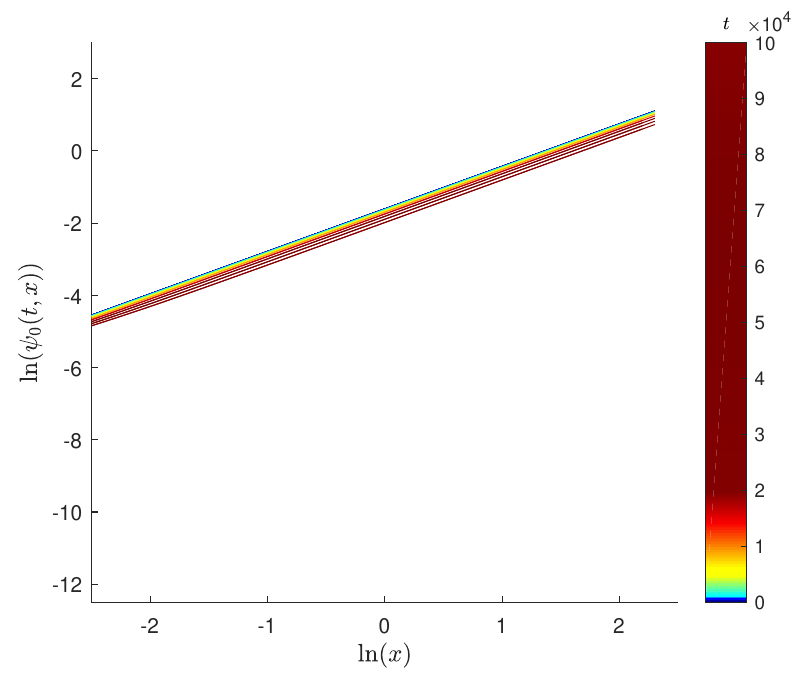}} 
	\end{minipage}
	{\begin{center}{\vspace{-1.1em} Subdiffusive case with $\mu = 0.3$ and $\beta (a) = \mu / (1+a)$. $\psi_0(0,x) = 0.2 (x-10)^{2/2-\mu}$. \vspace{-0.2em}}\end{center}}
	\caption{(Colour online) Decay of $\psi_0(t,\cdot)$ (left) and $\ln(\psi_0(t, \cdot))$ (right) for $\sigma = 1$, $t \in [0, 100000]$ (shown in the color bar) and $x \in [0, 20]$ with periodic boundary conditions.
		The presented plots are taken at $20$ regular intervals in $\ln(t) \in [0,11.5]$. As $t$ increases, each successive graph lies below the previous one for the larger values of $|x - 10|$.}
	\label{CGM_fig_decay_psi_0}
\end{figure}

For a visual representation of the evolution in time of the solution $\psi_0$ of the Hamilton-Jacobi equation~\eqref{CGM_eq_eq} in one space dimension, we refer to Figure~\ref{CGM_fig_decay_psi_0}, which is the result of a weighted essentially non-oscillatory (WENO) scheme of order 5 with Lax-Friedrichs numerical flux. We refer the reader to~\cite{Paul} for a review of such numerical methods. In Figure~\ref{CGM_fig_decay_psi_0}, the initial data taken for the first and second 
subfigures is the same, in order to illustrate how subdiffusion slows down significantly as time advances. 

The initial conditions in the first 
and third subfigures are chosen so as to decay with a preserved profile in 
$\rm{log}-\rm{log}$ scale for the Hamilton-Jacobi equations $\partial_t \psi + \tilde H(\partial_x \psi)$, with $\tilde H$ given by the approached expressions at $0$ of $H$~\eqref{CGM_eq_H_at_0}. Those are, respectively:
\begin{align}
&\partial_t \psi + |\partial_x \psi|^2 = 0, \qquad &\text{for the diffusive case,} \\
&\partial_t \psi + |\partial_x \psi|^{2/\mu} = 0, \qquad &\text{for the subdiffusive case}.\label{eq_heuristic_sd}
\end{align}
Injecting the Ansatz $\psi(t,x) = x^{\alpha}/t^\beta$ into the first equation yields 
$\ln (\psi) = 2 \ln x - \ln t$.
Injecting the same Ansatz into the second equation yields 
$\ln (\psi) = \frac 2 {2-\mu} \ln x - \frac \mu {2-\mu} \ln t$.

The values of $|\partial_x \psi|$ being low enough, the numerically computed solutions of the Hamilton-Jacobi equations exhibit a decay that agrees reasonably with our heuristic above.


	\section{Uniform local boundedness and Lipschitz continuity of $\psi_\e$}
	\label{CGM_sec_bounds}

\cred{For the sake of clarity, we present all our proofs in one-dimension of space $d=1$.
Extension to the higher dimensional case will be commented at crucial points throughout the proofs.}
	
	We will work over the set $[0,T] \times \RR$ for some $T > 0$, and we will denote by $C$ any positive real constant whose value is irrelevant.
%
	The subscript $\e$ may be dropped in the absence of confusion, for the sake of notations.

	This whole section deals with the proof of the following Theorem.

	\begin{thm} \label{CGM_thm_bounds}
		Let $T>0$ and $0 < \e < 1$. Under hypotheses \ref{CGM_hyp_omega} and \ref{CGM_hyp_IC}, $\psi_\e$ is bounded in $W^{1,\infty}_{\rm{loc}}([0,T] \times \RR)$ uniformly in $\e$, with the following quantitative bounds, where $(t,x) \in [0,T] \times \RR$:
		\begingroup
		\setlength{\abovedisplayskip}{0 pt plus 0 pt minus 0pt}
		\begin{enumerate}
			\item \label{CGM_thm_w1infty_psi_simple}
			\begin{align}
				\psi_\e (t,x) & \geq \inf v - \ln\beta(0) + \inf \eta , \label{CGM_eq_LBpsi_simple} \\
				\psi_\e (t,x) & \leq \inf(\phi_\e^0) + C_L |x| + C(\mu, \sigma, C_L, \| \eta \|_\infty) + (1+\mu)T. \label{CGM_eq_UBpsi_simple}
			\end{align}
			\item  \label{CGM_thm_w1infty_psi_x_simple}
			\begin{equation} \label{CGM_eq_Bpsi_x_simple} 
			\rm{Lip} ( \psi_\e (t, \cdot) )  \leq C_L.
			\end{equation}			
			\item  \label{CGM_thm_w1infty_psi_t_simple}
			\begin{equation} \label{CGM_eq_Bpsi_t_simple}
			{\rm Lip } \; \psi_\e (\cdot ,x) \leq \max \left( \mu ( 1 + \mu ) ,  \int_\RR \omega(z) \exp \left( C_L |z| \right) \,{\rm d}z \right) .
			\end{equation}
		\end{enumerate}
		\endgroup
	\end{thm}

	Subsection~\ref{CGM_ssec_psi} proves certain more accurate $\e$-dependent bounds~\eqref{CGM_eq_UB_psi_e_brut} from which the uniform bounds of Theorem~\ref{CGM_thm_bounds}.\ref{CGM_thm_w1infty_psi_simple} follow. The Lipschitz continuity results of Theorem~\ref{CGM_thm_bounds}.\ref{CGM_thm_w1infty_psi_x_simple} and Theorem~\ref{CGM_thm_bounds}.\ref{CGM_thm_w1infty_psi_t_simple} are proved in Subsections~\ref{CGM_ssec_psi_x} and~\ref{CGM_ssec_psi_t} respectively.
	
	\begin{rem}
		As mentioned previously in Subsection~\ref{CGM_ss_space_homo}, the space-homogeneous problem exhibits a self-similar decay in the original variables~\cite{BLM}. This precludes any time uniform $L^\infty$ bound of the solution $\psi_\e(t,x)$, as seen in the time correction in the upper bound \eqref{CGM_eq_UBpsi_simple}, which is more precisely of the form $\e (1 + \mu) \ln (1 + T/\e)$ as shown in the proof of the bound.
	\end{rem}

	\subsection{Local boundedness of $\psi_\e$} \label{CGM_ssec_psi}
	This subsection deals with the proof of Theorem~\ref{CGM_thm_bounds}.\ref{CGM_thm_w1infty_psi_simple}. 
	
	\begin{proof}[Proof of the lower bound~\eqref{CGM_eq_LBpsi_simple} of Theorem~\ref{CGM_thm_bounds}.\ref{CGM_thm_w1infty_psi_simple}]~
		
		From the scaling~\eqref{CGM_eq_def_phi} and the Ansatz~\eqref{CGM_eq_Ansatz} we see that $n^0$ can be expressed as follows:
		\[
		n^0(x,a) = n_\e^0(\e x, a) = \exp \left( -\frac {v(\e x)} \e - \eta(\e x, a) \right) \mathbb{1}_{[0,1]}(a),
		\]
		where we have dropped the subbscript $\e$.

		Let us define $\bar n : \RR_+ \times \RR_+ \to \RR_+$ as the solution of the following homogeneous problem:
		\begin{equation} \label{CGM_eq_aux_barn}
		\left\{
		\begin{aligned}
		&\partial_t \bar n (t,a) + \partial_a \bar n (t,a) + \beta(a) \bar n (t,a) = 0 \smallskip\\
		&\bar n (t,0) = \int_0^\infty \beta(a) \bar n (t,a) \,  \,{\rm d} a \smallskip\\
		&\bar n^0 (a) = \exp \left( - \inf_x v(\e x) / \e \right) \exp \left( - \inf_x \eta (\e x,a) \right) \mathbb{1}_{[0,1]}(a) \geq \sup_x n^0(x,a).
		\end{aligned} \right.
		\end{equation}
		Since $\beta$ is a non-increasing function, $\omega$ is a probability measure, $n$ is the solution of equation~\eqref{CGM_eq_n} for an initial condition $n^0$ and
		$\bar n^0 (a) \geq \sup_x n^0(x,a)$, it follows that for any $t,a \geq 0$:
		\[
		\bar n (t/\e,a) \geq \sup_x n(t/\e, x/\e, a) = \sup_x n_\e(t,x,a).
		\] 
		Moreover, since $\beta$ is non-increasing and the $L^1$ norm of $\bar n$ is preserved,
		\[
		\bar n (t/\e,0) = \int_0^\infty \beta(a) \bar n (t/\e,a) \,  \,{\rm d} a \leq \beta(0) \int_0^\infty \bar n (t/\e,a) \,  \,{\rm d} a = \beta(0) \|\bar n^0\|_{L^1}.
		\]
		It follows that:
		\[
		\psi_\e (t,x) = -\e \ln n_\e (t,x,0) \geq -\e \ln (\beta(0) \|\bar n^0\|_{L^1}).
		\]
		After computing the $L^1$ norm of $\bar n^0$:
		\[
		\|\bar n^0\|_{L^1} = \exp ( -\inf_x v / \e ) \int_0^\infty \exp (- \inf_x \eta(x,a) ) \,  \,{\rm d} a,
		\]
		and since $\eta$ is bounded (Hypothesis~\ref{CGM_hyp_IC}.\ref{CGM_hyp_IC_inf_eta}) and the integral over age is taken over $[0,1]$ due to the compactness of the initial support in age, we obtain the claimed result \eqref{CGM_eq_LBpsi_simple}. 
	\end{proof}

	The rest of this subsection is devoted to proving the upper bound~\eqref{CGM_eq_UBpsi_simple} in Theorem~\ref{CGM_thm_bounds}.\ref{CGM_thm_w1infty_psi_simple}.
	
	\begin{proof}[Proof of the upper bound~\eqref{CGM_eq_UBpsi_simple} of Theorem~\ref{CGM_thm_bounds}.\ref{CGM_thm_w1infty_psi_simple}]~

		From equation~\eqref{CGM_eq_eq} we recover:
		\begin{multline*}
		\exp \left( - \frac{1}{\e} \psi_\e (t,x) \right) = \int_0^{t/\e} \int_\RR \Phi(a) \omega(z) \exp \left(- \frac 1 \e \psi_\e(t-\e a, x - \e z) \right) \,{\rm d} z \,{\rm d} a \\
		+ \int_{0}^{1} \int_\RR \Phi \left(a + \frac t \e \right) \omega(z)  \exp\left(- \frac 1 \e \phi_\e^0(x - \e z, a ) \right) \exp \left( \int_0^{a} \beta \right) \,{\rm d} z \,{\rm d} a.
		\end{multline*}
		Since the first right-hand side term is non-negative and $v$ has at most linear growth, by Hypothesis~\ref{CGM_hyp_IC}.\ref{CGM_H2_sublinearity},
		\begin{align}\label{CGM_eq_UB_psi_e_brut_proof}
		e^{- \psi_\e (t,x) / \e} & \geq \int_0^1 \int_\RR \Phi \left(a + \frac{t}{\e}\right) \omega(z)  e^{- \frac 1 \e \left( \inf{\phi_\e^0} + C_L |x-\e z| \right)} e^{-\eta (x-\e z, a)} e^{\int_0^a \beta} \,{\rm d} z \,{\rm d} a \\
		& \geq \frac{\Phi \left(1 + \frac{t}{\e}\right)}{\Phi(0)} \int_0^1 \int_\RR \Phi(a)  e^{\int_0^a \beta} \omega(z) \exp \left( - \frac{\inf{\phi_\e^0}}{\e} - \frac{C_L |x|}{\e} - C_L |z| \right) \,{\rm d}z \,{\rm d}a \exp \left(- \|\eta\|_\infty \right) \nonumber \\
		& \geq e^{- \frac{\inf{\phi_\e^0}}{\e}} e^{- \frac{C_L |x|}{\e}} \frac{\Phi \left(1 + \frac{t}{\e}\right)}{\Phi(0)} \int_0^1  \Phi(a)  e^{\int_0^a \beta} \,{\rm d}a \int_\RR \frac{\exp \left( - \frac{z^2}{2 \sigma^2}  - C_L |z| \right)}{\sigma \sqrt{2 \pi}}  \,{\rm d}z \exp \left(- \|\eta\|_\infty \right). \nonumber
		\end{align}
		Moreover, the integral in age can be computed explicitly:
		\[
		\int_0^1 \Phi(a)  e^{\int_0^a \beta} \,{\rm d}a = \mu \ln(2),
		\]
		and the integral in space is bounded below by some constant $e^{-C}$, where $0 < C < \infty$ depends on $\sigma$ and $C_L$. Hence, we have the lower bound:
		\[
		e^{- \psi_\e (t,x) / \e} \geq  \exp \left(- \frac{\inf{\phi_\e^0}}{\e}\right) \exp \left(- \frac{C_L |x|}{\e}\right)  \frac{\Phi \left(1 + \frac{t}{\e}\right)}{\Phi(0)} \mu \ln(2)  \exp\left( - C \right)  \exp\left(- \|\eta\|_\infty \right).
		\]
		Taking the logarithm of the above expression yields an $\varepsilon$-dependent bound:
		\begin{equation}\label{CGM_eq_UB_psi_e_brut}
		\psi_\e (t,x) \leq \inf{\phi_\e^0} + C_L |x|  +  \e (1+\mu) \ln \left(1 + \frac{t}{\e}\right) + \e  C,
		\end{equation}
		for some positive $C$ depending on $\mu, \sigma, C_L, \| \eta \|_\infty$. Taking $0 < \e < 1$ leads in turn to the desired bound of equation~\eqref{CGM_eq_UBpsi_simple}.
	\end{proof}

	\subsection{Lipschitz continuity in $x$ of $\psi_\e$} \label{CGM_ssec_psi_x}

	This subsection deals with the proof of Theorem~\ref{CGM_thm_bounds}.\ref{CGM_thm_w1infty_psi_x_simple}.

	The keystone of our proof is an application of the maximum principle to the increase rate of $\psi_\e$. Let us set useful notations. Let $ h \in (0,1)$. We name the following differences:
	\begin{equation}\label{CGM_eq_def_psi_x_1}
	\left\{
	\begin{aligned}
	Z(t,x,a,z) = \psi_\e (t,x) - \psi_\e (t -\e a, x - \e z) \\
	Z_h(t,x,a,z) = Z(t,x+h,a,z) \\
	Y(t,x,a,z) = \psi_\e (t,x) - \phi_\e^0 (x - \e z, a - t / \e) \\
	Y_h(t,x,a,z) = Y(t,x+h,a,z).
	\end{aligned}
	\right.
	\end{equation}
	We define the difference quotients:
	\begin{equation}\label{CGM_eq_def_psi_x_2}
	\left\{
	\begin{aligned}
	u_h (t,x) = \frac{1}{h} \left(\psi_\e (t, x+h) - \psi_\e (t,x) \right), \\
	w_h^0(x,a) =  \frac{1}{h} \left(\phi_\e^0 (x+h, a) - \phi_\e^0 (x, a) \right).
	\end{aligned}
	\right.
	\end{equation}

	The use of the maximum principle requires bounded functions. As such, we introduce the following truncation (from above) of the initial data:
	\begin{equation}\label{CGM_eq_modified_IC}
	\left\{
	\begin{aligned}
	& \phi_\e^{0,R} (t,x) = v_\e^R (x) + \e \eta_\e (x,a) + \chi_{[0,1]}(a), \qquad 
	\\
	& v_\e^R(x) = \min \left( v_\e (x), \inf{\phi_\e^0} + C_L R \right).
	\end{aligned}
	\right.
	\end{equation}
	The  functions $\psi_\e^R$, $Z^R, Y^R, u_h^R$ and $w_h^{0,R}$ are defined accordingly.
	The upper bound in \eqref{CGM_eq_UBpsi_simple} becomes a true uniform bound, the term $C_L |x|$ being replaced with $ C_L R$. Additionaly, it is clear that, $\e$ being fixed, the original problem is recovered as $R\to +\infty$. Hence it is sufficient to prove the Lipschitz bound uniformly with respect to $R$.
	
	We begin with the proof of the upper bound. Assume by contradiction that there exist $\nu>0$ and $(t_0,x_0)$ such that~\footnote{meaning that $u_h^R(t_0,x_0)$ is greater than both terms.}
	\begin{equation}\label{eq:contradiction max pcple 1}
	u_h^R (t_0, x_0)    \geq 
	\begin{cases}
	\sup_{\RR \times [0,1]} w_h^{0,R} (x,a) + \nu\medskip\\
	\sup_{[0,T]\times \RR} u_h^R - \dfrac{\nu}{N}
	\end{cases}  \,,
	\end{equation}
	where $N$ is a constant depending on $T, h, \e, R, \nu$, to be determined below.

	By subtracting equation~\eqref{CGM_eq_eq}  evaluated at $(t,x)$ from the same equation evaluated at $(t, x+h)$ we get:
	\begin{multline*}
	0 = \int_0^{ \frac t \e} \int_\RR \Phi(a) \omega(z) \left[ \exp \left( Z_h^R(t,x,a,z) / \e \right) - \exp \left( Z^R(t,x,a,z) / \e \right) \right] \,{\rm d}z \,{\rm d}a \\
	+ \int_0^1 \int_\RR \Phi \left(a + \frac{t}{\e}\right) \omega(z) e^{\int_0^a \beta} \left[ \exp \left( Y_h^R(t,x,a,z) / \e \right) - \exp \left( Y^R(t,x,a,z) / \e \right) \right] \,{\rm d}z \,{\rm d}a.
	\end{multline*}
	By factoring out $Y$ and $Z$, we obtain the following identity at $(t_0,x_0)$:
	\begin{multline}\label{CGM_eq_discrete_max_pple_R}
	0 =  \int_0^{\frac {t_0} \e} \int_\RR \Phi(a) \omega(z) e^{ \frac{Z^R}{\e} } \left[\exp \left(\frac{h}{\e} \left[u_h^R(t_0,x_0) - u_h^R(t_0-\e a, x_0-\e z) \right]\right) - 1\right] \,{\rm d}z \,{\rm d}a \\
	+ \int_0^1 \int_\RR \Phi \left(a + \frac{t_0}{\e}\right) \omega(z) e^{\int_0^a \beta} e^{ \frac{Y^R}{\e}} \left[\exp \left(\frac{h}{\e} \left[u_h^R(t_0,x_0) - w_h^{0,R}(t_0-\e a, x_0-\e z) \right]\right) - 1\right] \,{\rm d}z \,{\rm d}a. \\
	\end{multline}
	Using \eqref{eq:contradiction max pcple 1}, we deduce that 
	\begin{align}
	0 \geq & \left ( \exp\left (- \dfrac{h}{\e} \dfrac{\nu}{N}\right ) - 1 \right )   \int_0^{\frac {t_0} \e} \int_\RR \Phi(a) \omega(z) e^{ \frac{Z^R}{\e} }  \,{\rm d}z \,{\rm d}a \nonumber\\
	& \quad +  \left ( \exp\left ( \dfrac{h}{\e} \nu\right ) - 1 \right )  \int_0^1 \int_\RR \Phi \left(a + \frac{t_0}{\e}\right) \omega(z) e^{\int_0^a \beta} e^{ \frac{Y^R}{\e}}  \,{\rm d}z \,{\rm d}a. \nonumber\\
	\geq & - \nu \left ( \dfrac{h}{\e N} \sup e^{ \frac{Z^R}{\e} }\right ) + \left ( \exp\left ( \dfrac{h}{\e} \nu\right ) - 1 \right ) \left ( \inf e^{ \frac{Y^R}{\e}} \right )\int_0^1 \int_\RR \Phi \left(a + \frac{T}{\e}\right) \omega(z) e^{\int_0^a \beta} \,{\rm d}z \,{\rm d}a.  \label{eq:contradict nu}
	\end{align}
	Therefore, it is possible to choose $N$ sufficiently large {\em  a priori}, (independently of $\nu\in (0,1)$), such that the right hand side in \eqref{eq:contradict nu}  is positive. This is a contradiction.

	This contradiction proves any space difference quotient of $\psi_\e$ is bounded above by the space Lipschitz constant of $\phi_\e^0$. The lower bound is proved in a similar way, and we recover Theorem~\ref{CGM_thm_bounds}.\ref{CGM_thm_w1infty_psi_x_simple}.

	\subsection{Lipschitz continuity in $t$ of $\psi_\e$} \label{CGM_ssec_psi_t}
	
	We proceed similarly for the time Lipschitz estimate. However, we bypass the rigorous use of difference quotients, as in \eqref{CGM_eq_def_psi_x_2}, but we differentiate the equation with respect to time. A rigorous proof can be obtained by a straightforward adaptation of the following arguments.
	
	We may reformulate \eqref{CGM_eq_eq} as follows, 
	\begin{multline} \label{CGM_eq_eq_sec_2}
	1 = \int_0^{t/\e} \Phi(a) \int_\RR \omega(z)  \exp\left( \frac 1 \e \left[ \psi^R_\e (t,x) - \psi^R_\e(t-\e a, x - \e z) \right] \right) \,{\rm d} z \,{\rm d} a \\
	+ \int_{0}^{1} \Phi(a+t/\e) \int_\RR \omega(z)  \exp\left(\frac 1 \e \left[ \psi^R_\e (t,x) - \phi_\e^{0,R}(x - \e z, a ) \right] + \int_0^{a} \beta \right) \,{\rm d} z \,{\rm d} a.
	\end{multline}
	and then we differentiate with respect to $t$ and multiply  by $\e e^{\frac 1 \e \psi_\e^R (t,x)}$ so as to get:
	\begin{equation} \label{CGM_eq_psi_t_details}
	\begin{aligned}
	0 & = \Phi \left( \frac t \e \right) 
	\int_\RR \omega(z) e^{-\frac 1 \e \psi^R(0, x-\e z)} \,{\rm d}z 
	\\
	& + \int_0^{t / \e} \int_\RR \Phi(a) \omega(z) \left[ \partial_t \psi_\e^R (t,x) - \partial_t \psi_\e^R (t-\e a, x-\e z) \right] e^{- \frac 1 \e \psi_\e^R (t - \e a, x - \e z)} \,{\rm d}z\,{\rm d}a \\
	& +  \int_{0}^{1} \int_\RR  \Phi\left (a+\frac t\e\right ) \omega(z) \left [ \partial_t \psi_\e^R (t,x) + \dfrac{\Phi'(a+\frac t\e)}{\Phi(a+\frac t\e)} \right ]  \exp\left( - \frac 1 \e  \phi_\e^{0,R}(x - \e z, a )  + \int_0^{a} \beta \right) \,{\rm d}z\,{\rm d}a.
	\end{aligned}
	\end{equation}

	We examine the upper bound and the lower bound separately. For the upper bound, assume by contradiction that there exists $\nu>0$, and $(t_0,x_0)$ such that
	\begin{equation}\label{eq:contradiction max pcple 2}
	\partial_t\psi_\e^R (t_0, x_0)    \geq 
	\begin{cases}
	\mu(1 + \mu) + \nu\medskip\\
	\sup_{[0,T]\times \RR} \partial_t\psi_\e^R  - \dfrac{\nu}{N}
	\end{cases}  \,,
	\end{equation}
	where $N$ is a constant depending on $T, h, \e, R, \nu$, to be determined below.~\footnote{Equations~\eqref{eq:contradiction max pcple 2} and~\eqref{eq:contradiction max pcple 3} assume $\partial_t\psi_\e^R$ is bounded. This is true for the corresponding difference quotient at any fixed time step. Difference quotients also simplify the generalisation to $d>1$.} By ignoring the first (positive) contribution, we deduce from \eqref{CGM_eq_psi_t_details} that
	\begin{equation} \label{CGM_eq_psi_t_details2}
	\begin{aligned}
	0 & \geq 
	-\dfrac{\nu}{N} \left ( \sup e^{- \frac 1 \e \psi_\e^R }  \right )  \\
	& + \nu \left ( \inf \exp\left(- \frac 1 \e  \phi_\e^{0,R}  \right)  \right ) \int_{0}^{1} \int_\RR  \Phi\left (a+\frac {T}\e\right ) \omega(z)   e^{\int_0^a \beta}   \,{\rm d}z\,{\rm d}a.
	\end{aligned}
	\end{equation}
	where we have used the fact that 
	\begin{equation}
	\frac{\Phi' \left( a + \frac {t_0} \e \right)}{\Phi \left( a + \frac {t_0} \e \right)} = \frac{\mu (1+\mu)}{1 + a + \frac {t_{0}} \e} \leq \mu(1+\mu)\, .
	\end{equation}
	Again, by choosing the constant $N$ large enough, we reach a contradiction.

	For the lower bound, we can ignore the contribution involving $\frac{\Phi'}{\Phi}$ as it is negative. We assume by contradiction that there exists $\nu>0$, and $(t_0,x_0)$ such that
	\begin{equation}\label{eq:contradiction max pcple 3}
	\partial_t\psi_\e^R (t_0, x_0)    \leq 
	\begin{cases}
	- \int_\RR \omega(z) \exp \left( C_L |z| \right) \,{\rm d}z  -  \nu\medskip\\
	\inf_{[0,T]\times \RR} \partial_t\psi_\e^R  + \dfrac{\nu}{N}
	\end{cases}  \,,
	\end{equation}
	where $N$ is a constant depending on $T, h, \e, R, \nu$, to be determined below. The trick here is to make appear comparable quantities in \eqref{CGM_eq_psi_t_details}:
	\begin{equation} \label{eq:cancellation}
	\begin{aligned}
	0 & \leq  \Phi \left( \frac {t_0} \e \right) 
	\int_\RR \omega(z) \left ( \int_{0}^{1} \int_\RR  \Phi(a)\omega(y) e^{-\frac 1 \e \phi_\e^{0,R}(x-\e z -\e y,a)} \, {\rm d}y\, {\rm d}a \right ) \,{\rm d}z 
	+ \dfrac\nu N \left ( \sup e^{- \frac 1 \e \psi_\e^R }  \right ) \\
	& + \left [ - \int_\RR \omega(z) \exp \left( C_L |z| \right) \,{\rm d}z  - \nu  \right ] \int_{0}^{1} \int_\RR  \Phi\left (a+\frac {t_0}\e\right ) \omega(z)   \exp\left( - \frac 1 \e  \phi_\e^{0,R}(x - \e z, a )  + \int_0^{a} \beta \right) \,{\rm d}z\,{\rm d}a\\
	&\leq \left ( \int_\RR \omega(y) \exp \left( C_L |y| \right) \,{\rm d}y \right )  \int_0^1 \int_\RR \Phi\left (a+\frac {t_0}\e\right ) \omega(z) \exp\left( - \frac 1 \e  \phi_\e^{0,R}(x - \e z, a )  + \int_0^{a} \beta \right) \,{\rm d}z\,{\rm d}a \\
	& + \dfrac\nu N \left ( \sup e^{- \frac 1 \e \psi_\e^R }  \right ) +  \left [ - \int_\RR \omega(z) \exp \left( C_L |z| \right) \,{\rm d}z  - \nu  \right ] \int_{0}^{1} \int_\RR  \Phi\left (a+\frac {t_0}\e\right ) \omega(z) \cdot \\
	& \quad \cdot \exp\left( - \frac 1 \e  \phi_\e^{0,R}(x - \e z, a )  + \int_0^{a} \beta \right) \,{\rm d}z\,{\rm d}a 
	\end{aligned}
	\end{equation}
	where we have used the following pointwise inequality which holds for any $t_0\geq 0$:
	\begin{equation}
	\Phi \left( \frac {t_0} \e \right)  \Phi(a) \leq \Phi \left( a +  \frac {t_0} \e \right)  \exp\left(   \int_0^{a} \beta \right)\, .
	\end{equation}
	By choosing $N$ sufficiently large, we arrive to a contradiction due to cancellations in \eqref{eq:cancellation}.

	\section{Viscosity limit procedure}\label{CGM_sec_viscosity_solution}
	
	In this section, we continue to work over $[0,T] \times \RR$. 
	
	We deduce from the Lipschitz estimates  that there exists a Lipschitz function $\psi_0$ such that $\psi_{\e} \to\psi_0$ locally uniformly, up to extraction. 
	We shall prove that $\psi_0$ is the unique viscosity solution of the Hamilton-Jacobi equation~\eqref{CGM_eq_limiting_HJ}, which we recall here:
	\[
	1 = \int_0^\infty \Phi(a) \exp\left(a \partial_t \psi_0 (t,x)\right) \,{\rm d}a \int_\RR \omega(z) \exp\left(z \partial_x \psi_0 (t,x)\right) \,{\rm d} z
	\]
	with initial condition $\psi_0(0,x) = v(x)$.

	Equation \eqref{CGM_eq_eq} is equivalent to the following, which allows us to define $\mathfrak A_\e$ and $\mathfrak B_\e$ and is better suited for the following proofs:
	\begin{equation}\label{CGM_eq_def_A_e}
	\begin{aligned}
	1 & = \left(\mathfrak{A}_{\e} + \mathfrak{B}_\e\right) (\psi_\e)(t,x), \qquad \text{where:} \smallskip\\
	& \mathfrak{A}_\e (\psi_\e)(t,x) = \int_0^{t / \e} \int_\RR \omega(z) \Phi(a) \exp \left( \frac 1 \e \left[\psi_\e(t,x) - \psi_\e(t-\e a, x - \e z)\right] \right) \,{\rm d}z\,{\rm d}a, \smallskip\\
	& \mathfrak{B}_\e (\psi_\e)(t,x) = \int_0^1 \int_\RR \omega(z)  \Phi(a+t/\e)  \exp \left( \frac 1 \e \left[\psi_\e(t,x) - \phi_\e^0(x - \e z, a) \right] \right) \exp \left( \int_0^a \beta \right) \,{\rm d}z\,{\rm d}a.
	\end{aligned}
	\end{equation}

	\subsection{Viscosity subsolution}
	
	\begin{prop}\label{CGM_prop_viscosity_subsolution}
		Under hypotheses~\ref{CGM_hyp_omega} and~\ref{CGM_hyp_IC}, $\psi_0$ is a viscosity subsolution of~\eqref{CGM_eq_limiting_HJ}. 
	\end{prop}
	
	\begin{proof}
		Let $\Psi \in \mathcal C^2(\RR^+ \times \RR)$ be a test function such that $\psi_0 - \Psi$ admits a maximum at $(t_0, x_0)$, with $t_0 >0$.
		By compactness in $W^{1,\infty}_{\rm{loc}}([0,T] \times \RR)$, thanks to the a priori estimates, we obtain for a subsequence of $\varepsilon \to 0$ which we will not rename: 
		$(t_\varepsilon, x_\varepsilon) \underset{\varepsilon \to 0}{\longrightarrow} (t_0, x_0)$, where $(t_\varepsilon, x_\varepsilon)$ is a point at which $\psi_\e - \Psi$ reaches its maximum. We have then:\\
		$\forall \varepsilon > 0 \; \forall (z,a) \in \RR \times [0, \frac{t_\varepsilon}{\varepsilon}]$, 
		\[
		\psi_\varepsilon (t_\varepsilon, x_\varepsilon) - \Psi (t_\varepsilon, x_\varepsilon) \geq \psi_\varepsilon (t_\varepsilon - \varepsilon a, x_\varepsilon - \varepsilon z) - \Psi (t_\varepsilon - \varepsilon a, x_\varepsilon - \varepsilon z).
		\]
		Since $\mathfrak{B}_\e$ is non-negative, it follows that:
		\[ 1 \geq \mathfrak A_\varepsilon (\psi_\varepsilon) (t_\varepsilon, x_\varepsilon) \geq \mathfrak A_\varepsilon (\Psi) (t_\varepsilon, x_\varepsilon) .\]
		However:
		\begin{equation} \label{CGM_eq_Taylor_expansion_Psi}
		\begin{aligned}
		\Psi (t_\varepsilon, x_\varepsilon) & - \Psi (t_\varepsilon - \varepsilon a, x_\varepsilon - \varepsilon z) =
		\varepsilon a \partial_t \Psi (t_\varepsilon, x_\varepsilon) + \varepsilon z \partial_x \Psi (t_\varepsilon, x_\varepsilon) \\
		& + \frac{1}{2}\varepsilon^2 \int_0^1 (1-s)^2 \left[ a^2 \partial_t ^2 \Psi (t_\varepsilon - \varepsilon s a , x_\varepsilon - \varepsilon s z) + 2 a z \cdot \partial_t \partial_x \Psi (t_\varepsilon - \varepsilon s a , x_\varepsilon - \varepsilon s z) \right. \\ 
		& \left. + z^2 \partial_x ^2 \Psi (t_\varepsilon - \varepsilon s a , x_\varepsilon - \varepsilon s z) \right] \,{\rm d}s .
		\end{aligned}
		\end{equation}
		Therefore we have, for all $A > 0$:
		\begin{align*}
		1 & \geq \int_0^A \int_{-A}^A \Phi(a) \omega(z) \exp \bigg\{ 
		a \partial_t \Psi (t_\varepsilon, x_\varepsilon) + z \partial_x \Psi (t_\varepsilon, x_\varepsilon)  \\
		& \left. + \frac{1}{2}\varepsilon \int_0^1 (1-s)^2     \left[     a^2 \partial_t ^2 \Psi (t_\varepsilon - \varepsilon s a , x_\varepsilon - \varepsilon s z) 
		+ 2 a z \cdot \partial_t \partial_x \Psi (t_\varepsilon - \varepsilon s a , x_\varepsilon - \varepsilon s z) \right. \right. \\
		& \left. +
		z^2 \partial_x^2 \Psi (t_\varepsilon - \varepsilon s a , x_\varepsilon - \varepsilon s z)     \right]     \,{\rm d}s \bigg\} \,{\rm d}z \,{\rm d}a .
		\end{align*}
		
		Since $\Psi$ is $\mathcal C^2$, the previous expression tends, for fixed $A$, when $\varepsilon \to 0$, to:
		
		\[
		1 \geq \int_0^A \int_{-A}^A \Phi(a) \omega(z) \exp \left[ a \partial_t \Psi (t_0, x_0) + z \partial_x \Psi (t_0, x_0) \right] \,{\rm d}z \,{\rm d}a .
		\]
		
		It follows that:
		
		\[
		1 \geq \int_0^\infty \int_\RR \Phi(a) \omega(z) \exp \left[ a \partial_t \Psi (t_0, x_0) + z \partial_x \Psi (t_0, x_0) \right] \,{\rm d}z \,{\rm d}a .
		\]
		
		Therefore $\psi_0$ is a viscosity subsolution of~\eqref{CGM_eq_limiting_HJ}. 
	\end{proof}

	\subsection{Viscosity supersolution} \label{CGM_ssec_visc_supersol}

	In order to prove that $\psi_0$ is a viscosity supersolution of~\eqref{CGM_eq_limiting_HJ}, 
	we need to control the $\mathfrak B_\e$ term in equation~\eqref{CGM_eq_def_A_e}, whose positivity sufficed in the previous subsection. This is tantamount to controlling the fate of the aging particles that come from the initial data and have never jumped.
	
	We proceed in several steps. The key idea is to compare the relative weigths of $\mathfrak{A}_\e$ and $\mathfrak{B}_\e$, by means of the quantity $\psi_\e(t,x) - \psi_\e(0,x)$. Because the sum of the two contributions equals one, we shall deduce that $\mathfrak{A}_\e\to 1$, and $\mathfrak{B}_\e\to 0$. Interestingly enough, we get a quantitative estimate on the convergence rate.
	
	\paragraph*{Step 1: A crude estimate on $\mathfrak{B}_\e$}
	The following Lemma boils down the estimate on $\mathfrak{B}_\e$ to some estimate on time increments of $\psi_\e$. Here, the boundedness of the age support is crucial.

	\begin{lem}[Simple bounds for $\mathfrak B_\e$] \label{CGM_lem_B_e}
		\begin{equation}
		\frac{\Phi(t/\e)}{\Phi(0)} \exp\left(\frac 1 \e \left[ \psi_\e(t,x) - \psi_\e (0,x) \right] \right) \leq \mathfrak B_\e(\psi_\e)(t,x) \leq \frac{\Phi(1+t/\e)}{\Phi(1)} \exp\left(\frac 1 \e \left[ \psi_\e(t,x) - \psi_\e (0,x) \right] \right) 
		\end{equation}
	\end{lem}
	\begin{proof}
		This is a consequence of the following claim: for all $h > 0$, 
		\[ a \mapsto \Phi(a+h) / \Phi(a) \] is an increasing function.
		Indeed,
		\begin{align*}
		& \frac{\rm d}{{\rm d} a} \frac{\Phi(a+h)}{\Phi(a)} = \frac{\Phi'(a+h) \Phi(a) - \Phi'(a)\Phi(a+h)}{(\Phi(a))^2} \\
		= & \frac{\exp\left(\int_0^a \beta \right) \exp\left(\int_0^{a+h} \beta \right)}{(\Phi(a))^2} \Big( \beta'(a+h)\beta(a) - \beta'(a)\beta(a+h) + \underbrace{\beta(a)\beta(a+h)[\beta(a) - \beta(a+h)]}_{\geq 0} \Big),
		\end{align*}
		which is positive since, $\beta(a) =\mu / (1+a)$ being non-increasing and convex, $\beta'(a+h) \geq \beta'(a)$ by convexity and $\beta(a) \geq \beta(a+h) \geq 0$. This proves the claim.
		
		We now write $\mathfrak B_\e$ as follows:
		\[
		\mathfrak B_\e(\psi_\e)(t,x)  =  \int_0^1  \frac{\Phi(a+t/\e)}{\Phi(a)} \Phi(a) \int_\RR \omega(z) \exp\left( \frac 1 \e \left[ \psi_\e(t,x) - \phi_\e^0 (x - \e z, a) \right] \right) \exp\left(\int_0^a \beta \right) \,{\rm d}z\,{\rm d}a
		\]
		and recover the lower and upper bounds by monotonicity and thanks to~\eqref{CGM_eq_eq}.
	\end{proof}
	
	\paragraph*{Step 2: A lower bound for $\mathfrak{A}_\e$}
	The goal of the following Lemma is to remove the $x$ variations from the contribution in $\mathfrak{A}_\e$. Hence, the problem will be reduced to estimate for a given $x$. This strongly relies on semi-concavity.  
	
	Semi-concavity is a natural regularity for Hamilton-Jacobi equations. It can result either from the propagation of regularity on the initial data, or on regularization property of the Hamilton-Jacobi equation~\cite[Chapter~3.3]{Evans}. The latter usually relies on uniform convexity of the Hamiltonian, which is not the case here. Below, we derive propagation estimates for $\e>0$.

	\begin{lem}[Lower bound for $\mathfrak A_\e$] \label{CGM_lem_A_e}
		For $\e$ small enough and $(t,x) \in [0,T] \times \RR$,
		\begin{equation}\label{CGM_eq_LB_A_e}
		\mathfrak{A}_\e(\psi)(t,x) 
		\geq \Bigl[1-\e\frac {\mathfrak C_{xx}} 2\int_\RR\omega(z)z^2\,{\rm d}z\Bigr]\int_0^{t/\e} \Phi(a) \exp \left( \frac 1 \e \left[ \psi_\e (t,x) - \psi_\e (t- \e a, x)\right] \right)\,{\rm d}a
		\end{equation}
		where $\mathfrak C_{xx}$ is the upper bound of $\partial_x^2 \phi_\e^0$ from Hypothesis~\ref{CGM_hyp_IC}.\ref{CGM_hyp_IC_semi_concavity}.
	\end{lem}
	\begin{proof}
		First, let us prove that the semi-concavity of the initial condition is preserved. By differentiating \eqref{CGM_eq_eq} twice with respect to $x$, we obtain:
		\begin{equation}
		\begin{aligned}
		0 =
		& \int_0^{t / \e} \int_\RR \omega(z) \Phi(a) \bigg[ \left( \partial_x \psi_\e(t,x) - \partial_x \psi_\e (t - \e a, x - \e z) \right)^2 + \\
		& \frac 1 \e \left(  \partial_x^2 \psi_\e(t,x) -  \partial_x^2 \psi_\e (t - \e a, x - \e z) \right) \bigg]
		\exp\left( \frac 1 \e \left[ \psi_\e(t,x) - \psi_\e (t - \e a, x - \e z) \right] \right) \,{\rm d}z\,{\rm d}a \medskip\\
		& + \int_0^1 \int_\RR \omega(z) \Phi ( a + t / \e ) \bigg[ \left( \partial_x \psi_\e(t,x) - \partial_x \phi_\e^0 (x - \e z , a) \right)^2 +\\ 
		&  \frac 1 \e \left(  \partial_x^2 \psi_\e(t,x) -  \partial_x^2 \phi_\e^0 (x - \e z , a) \right) \bigg] \exp\left(\frac 1 \e \left[\psi_\e(t,x) - \phi_\e^0 (x - \e z , a)\right]\right) \exp\left(\int_0^a \beta\right) \,{\rm d}z\,{\rm d}a,
		\end{aligned}
		\end{equation}
		Since $\psi_\e$ and $\phi_\e^0$ are Lipschitz continuous in $x$ and thanks to Rademacher's theorem they are almost everywhere differentiable, the squared terms are well defined and non-negative. Moreover, Hypothesis~\ref{CGM_hyp_IC}.\ref{CGM_hyp_IC_semi_concavity} gives us $\partial_x^2 \phi_\e^0 \leq \mathfrak C_{xx}$ in the sense of distributions. We recover an upper bound for $\partial_x^2 \psi_\e$. Indeed,  at $(t_0, x_0) = \arg \max  \partial_x^2 \psi_\e$~\footnote{\label{CGM_f_max} $\partial_x^2 \psi_\e$ may not reach its maximum but in this case it suffices to proceed as in subsection \ref{CGM_ssec_psi_x}.}, an application of the maximum principle allows us to recover:
		\begin{equation}\label{CGM_eq_dxx_psi}
		\partial_x^2 \psi_\e (t_0, x_0) \leq \mathfrak C_{xx}.
		\end{equation}
		Secondly, we deduce the following simple Taylor estimate, 
		\begin{align*}
		\exp \left(- \frac 1 \e \left[ \psi_\e ( t - \e a, x - \e z ) - \psi_\e (t - \e a, x) \right]\right)&\geq1-\frac 1 \e \left[ \psi_\e ( t - \e a, x - \e z ) - \psi_\e (t - \e a, x) \right]\\
		&\geq1-\frac 1 \e \left[ -\e z\,\partial_x\psi_\e ( t - \e a, x ) + \frac {\mathfrak C_{xx}} 2 \e^2z^2 \right].
		\end{align*}
		Then, since $\int_\RR z\omega(z)\,{\rm d}z=0$, we have,
		\begin{equation}\label{CGM_eq_DL_psi_loc}
		\int_\RR \omega(z) \exp\left(- \frac 1 \e \left[ \psi_\e ( t - \e a, x - \e z ) - \psi_\e (t - \e a, x) \right]\right)\,{\rm d}z\geq1-\e\frac {\mathfrak C_{xx}} 2\int_\RR \omega(z)z^2\,{\rm d}z.
		\end{equation}
		The result of the Lemma follows.
	\end{proof}
	
	\paragraph*{Step 3: An upper bound on $\psi_\e(t,x) - \psi_\e(0,x)$}
	
	We are now ready to apply the maximum principle on the time increment for a fixed $x$. 
	
	\begin{lem}[Upper bound] \label{CGM_lem_m_e}
		Let us fix $x \in \RR$. Let $m$ be the maximum over $t \in [0,T]$ of $\psi_\e ( t ,x) - \psi_\e(0,x)$. For $K =  (\mathfrak C_{xx} / 2) \int_\RR\omega(z)z^2\,{\rm d}z$, we have:
		\begin{equation}
		e^{m/\e} \leq 1 + \frac T \e + K \e \left(1 + \frac T \e \right)^{1+\mu}.
		\end{equation}
	\end{lem}
	\begin{proof} 
		We deduce from the identity $\mathfrak{A}_\e+\mathfrak{B}_\e=1$, from Lemma~\ref{CGM_lem_A_e} and from the lower bound in Lemma~\ref{CGM_lem_B_e}, that,
		\[
		1 \geq [1- K \e ] \int_0^{t/\e} \Phi(a) \exp \left( \frac 1 \e \left[ \psi_\e (t,x) - \psi_\e (t- \e a, x)\right] \right)\,{\rm d}a
		+ \frac{\Phi(t/\e)}{\Phi(0)} \exp\left(\frac 1 \e \left[ \psi_\e(t,x) - \psi_\e (0,x) \right] \right).
		\]
		Applying the maximum principle and denoting $t_0 = \arg \max_{[0,T]} \psi_\e (\cdot, x)$ results in
		\[
		1 \geq [ 1 - K \e ] \int_0^{t_0/\e} \Phi(a)\,{\rm d}a+ \frac{\Phi(t_0/\e)}{\Phi(0)} e^{m/\e},
		\]
		hence:
		\begin{align*}
		e^{m/\e} & \leq \frac{\Phi(0)}{\Phi(t_0 / \e)} \left[ 1 - (1-K\e) \int_0^{t_0/\e} \Phi(a) \, {\rm d}a \right] \\
		& \leq \left( 1 + \frac{t_0} \e \right)^{1+\mu} \left[ \int_{t_0/\e}^\infty \Phi(a) \,{\rm d}a + K \e \int_0^{t_0 / \e} \Phi(a) \,{\rm d}a \right] \\
		& \leq 1 + \frac { t_0} \e + K \e \left[\left(1 + \frac{t_0} \e \right)^{1+\mu} - \left(1 + \frac{t_0} \e \right)\right],
		\end{align*}
		hence the result.
	\end{proof}
	
	\paragraph*{Step 4: An upper bound on $\mathfrak{B}_\e$}
	
	Back to the upper bound in Lemma~\ref{CGM_lem_B_e}, we are in position to conclude.

	\begin{prop}[Upper bound for $\mathfrak B_\e$] \label{CGM_prop_B_e_to_0}
		Under hypotheses \ref{CGM_hyp_omega} and \ref{CGM_hyp_IC}, for any $(t,x) \in [0,T] \times \RR$, $\mathfrak B_\e$ decays in the following way as $\e \to 0$:
		\begin{equation} \label{CGM_eq_B_e}
		\mathfrak B_\e  \leq \e^\mu \frac{2^{1+\mu}}{t^{1+\mu}} \left[T + K_1 \e^{1-\mu} + K_2 \e\right]\,,
		\end{equation}
		for some explicit constants $K_1,K_2$.
	\end{prop}
	
	\begin{proof}
		Lemma~\ref{CGM_lem_m_e} and the upper bound in Lemma~\ref{CGM_lem_B_e} give us:
		\begin{align*}
		\mathfrak B_\e(\psi_\e)(t,x) & \leq \frac{\Phi(1+t/\e)}{\Phi(1)} \exp\left(\frac 1 \e \left[ \psi_\e(t,x) - \psi_\e (0,x) \right] \right) \\
		& \leq \frac{2^{1+\mu}}{(2+t/\e)^{1+\mu}} e^{m/\e} \\
		& \leq 2^{1+\mu} \left(\dfrac{\e}{t}\right )^{1+\mu}  \left( 1 + \frac T \e + K \e \left(1 + \frac T \e \right)^{1+\mu} \right)\,.
		\end{align*}
	\end{proof}

	\paragraph*{Step 5: Conclusion of the proof}
	
	The accurate upper bound on $\mathfrak{B}_\e$ that we have just proved allows us to proceed to the crucial result of this section.
	\begin{prop}\label{CGM_prop_viscosity_supersolution}
		Under hypotheses \ref{CGM_hyp_omega} and \ref{CGM_hyp_IC}, $\psi_0$ is a viscosity supersolution of~\eqref{CGM_eq_limiting_HJ}. 
	\end{prop}
	
	\begin{proof}
		Let $\Psi \in \mathcal C^2(\RR^+ \times \RR)$ be a test function such that $\psi_0 - \Psi$ admits a strict local minimum at $(t_0, x_0)$, with $t_0 >0$.  
		We make the distinction between two cases:\medskip\\
		
		\subparagraph*{If $\partial_t\Psi(t_0,x_0)\geq 0$,} then we get immediately 
		\[
		1 \leq \int_0^\infty \Phi (a) \exp \left( a \partial_t \Psi(t_0,x_0) \right)  \,{\rm d}a  \int_\RR  \omega(z)  \exp \left( z \partial_x \Psi (t_0, x_0) \right) \,{\rm d}z.
		\]
		Indeed, if $\partial_t\Psi(t_0,x_0)> 0$ then the right hand side is infinite. Whereas, if $\partial_t\Psi(t_0,x_0)=0$, this equality follows from the symmetry of $\omega$, see also \eqref{CGM_eq_aux_omega}. \medskip\\
		
		\subparagraph*{If $\partial_t\Psi(t_0,x_0)< 0$,} then there exists $\nu>0$, and a ball of radius $0<2h<t_0/10$, $B((t_0,x_0),2h)$ such that $\partial_t\Psi(t,x)<-\nu$ over the ball. 
		On the other hand, by uniform convergence of $\psi_\e$ to $\psi_0$, there exists $(t_\varepsilon, x_\varepsilon)$ such that $\psi_\e - \Psi$ reaches a local minimum at $(t_\varepsilon, x_\varepsilon)$. We assume that $\e$ is small enough such that $(t_\varepsilon, x_\varepsilon)\in B((t_0,x_0),h)$.
		
		The contribution $\mathfrak B_\varepsilon$ is handled thanks to Proposition~\ref{CGM_prop_B_e_to_0}, uniformly in $t \in [t_0 / 2, T]$: 
		\[
		\forall \delta > 0, \; \exists \, \varepsilon_\delta > 0 \;  | \; \forall \varepsilon \in (0,\varepsilon_\delta),  \; \mathfrak B_\varepsilon (\Psi)(t_\varepsilon, x_\varepsilon) < \delta .
		\]
		The contribution $\mathfrak{A}_\e$ is handled by splitting the time integral into two contributions: those ages which are smaller than $h/\e$, and those ages which are greater. The small ages are dealt with thanks to the local minimum property: 
		\begin{equation*} \forall (a,z) \in B((t_\e,x_\e),h/\e)\quad
		\psi_\varepsilon (t_\varepsilon, x_\varepsilon) - \Psi (t_\varepsilon, x_\varepsilon) \leq \psi_\varepsilon (t_\varepsilon - \varepsilon a, x_\varepsilon - \varepsilon z) - \Psi (t_\varepsilon - \varepsilon a, x_\varepsilon - \varepsilon z).
		\end{equation*}
		Recalling the identity $\mathfrak{A}_\e+\mathfrak{B}_\e =1$, we deduce
		\begin{equation}\label{CGM_eq_I_II}
		1 - \delta \leq  A_\varepsilon (\psi_\varepsilon) (t_\varepsilon, x_\varepsilon) 
		=  I + II + III, 
		\end{equation}
		where we set  $h>0$ and define:
		\begin{equation}
		\begin{aligned}
		& I = \int_0^{\frac h  \e} \int_{ - \frac h \e}^{\frac h \e } \Phi(a) \omega(z) \exp \left( \frac 1 \e \left[\psi_\e(t_\e,x_\e) - \psi_\e(t_\e-\e a, x_\e - \e z)\right] \right) \,{\rm d}z\,{\rm d}a,\\
		& II = \int_0^{\frac h \e} \int_{\RR \setminus \left[- \frac h \e , \frac h \e \right]} \Phi(a) \omega(z) \exp \left( \frac 1 \e \left[\psi_\e(t_\e,x_\e) - \psi_\e(t_\e-\e a, x_\e - \e z)\right] \right) \,{\rm d}z\,{\rm d}a,\\
		& III = \int_{\frac h \e}^{1 + \frac t \e} \int_\RR \Phi(a) \omega(z) \exp \left( \frac 1 \e \left[\psi_\e(t_\e,x_\e) - \psi_\e(t_\e-\e a, x_\e - \e z)\right] \right) \,{\rm d}z\,{\rm d}a.
		\end{aligned}
		\end{equation}

		$\bullet$ Limit of $I$ -- small ages and spaces.
		
		Thanks to the local maximum property we have:
		\[
		I \leq \int_0^{\frac h  \e} \int_{ - \frac h \e}^{\frac h \e } \Phi(a) \omega(z) \exp \left( \frac 1 \e \left[\Psi(t_\e,x_\e) - \Psi(t_\e-\e a, x_\e - \e z)\right] \right) \,{\rm d}z\,{\rm d}a,
		\]
		on which we perform the same Taylor expansion as in~\eqref{CGM_eq_Taylor_expansion_Psi}, which yields:
		\begin{align*}
		I & \leq \int_0^{\frac h \e} \int_{- \frac h \e}^{\frac h \e} \Phi(a) \omega(z) \exp \bigg\{ 
		a \partial_t \Psi (t_\varepsilon, x_\varepsilon) + z \partial_x \Psi (t_\varepsilon, x_\varepsilon)  \\
		& \left. + \frac{1}{2}\varepsilon \int_0^1 (1-s)^2     \left[     a^2 \partial_t ^2 \Psi (t_\varepsilon - \varepsilon s a , x_\varepsilon - \varepsilon s z) 
		+ 2 a z \cdot \partial_t \partial_x \Psi (t_\varepsilon - \varepsilon s a , x_\varepsilon - \varepsilon s z) \right. \right. \\
		& \left. + z^2 \partial_x^2 \Psi (t_\varepsilon - \varepsilon s a , x_\varepsilon - \varepsilon s z)     \right]     \,{\rm d}s \bigg\} \,{\rm d}z \,{\rm d}a .
		\end{align*}
		Since $\Psi \in \mathcal{C}^2$ only takes values over $B( (t_0, x_0), 2 h)$ in the expression above, uniformly in $\e$, a domination argument allows us to pass to the limit $\e \to 0$ and recover the following limit for the right hand side:
		\[
		\int_0^\infty \int_\RR \Phi (a) \omega(z) \exp \left( [a \partial_t \Psi + z \partial_x \Psi](t_0, x_0)\right) \,{\rm d}z \,{\rm d}a.
		\]

		$\bullet$ Limit of $II$ -- small ages, large spaces.
		
		Since $\psi_\e$ is Lipschitz continuous in $x$ with some constant $L$, we can localise the expression of $II$ at $x_\e$ at a price:
		\[
		II \leq  \int_0^{h / \e} \int_{\RR \setminus \left[- \frac h \e , \frac h \e \right]} \Phi(a) \omega(z) e^{L |z|} \exp \left( \frac 1 \e \left[\psi_\e(t_\e,x_\e) - \psi_\e(t_\e-\e a, x_\e)\right] \right) \,{\rm d}z\,{\rm d}a.
		\]
		Thanks to the local maximum property,
		\[
		II \leq  \int_0^{h / \e} \int_{\RR \setminus \left[- \frac h \e , \frac h \e \right]} \Phi(a) \omega(z) e^{L |z|} \exp \left( \frac 1 \e \left[\Psi(t_\e,x_\e) - \Psi(t_\e-\e a, x_\e)\right] \right) \,{\rm d}z\,{\rm d}a.
		\]
		And by negativity of $\partial_t \Psi$ around $(t_0, x_0)$,
		\[
		II \leq  \int_0^{h / \e} \int_{\RR \setminus \left[- \frac h \e , \frac h \e \right]} \Phi(a) \omega(z) e^{L |z|} e^{-\nu a} \,{\rm d}z\,{\rm d}a,
		\]
		which converges to $0$ as $\e \to 0$.
		
		$\bullet$ Limit of $III$ -- large ages.
		
		Since $\psi_\e$ is Lipschitz continuous with some Lipschitz constant $L$, we recover:
		\[
		III  \leq \int_{\frac h \e}^{1 + \frac t \e} \int_\RR \Phi \left(a \right) \omega(z) e^{|z|L} \exp \left( \frac 1 \e \left[\psi_\e(t_\e,x_\e) - \psi_\e(t_\e -\e a, x_\e )\right] \right) \,{\rm d}z\,{\rm d}a.
		\]
		We have:
		\begin{multline*}
		\psi_\e(t_\e,x_\e) - \psi_\e(t_\e -\e a, x_\e ) = \\
		\psi_\e(t_\e,x_\e) -  \psi_\e(t_\e - h, x_\e) + \Psi(t_\e - h, x_\e) - \Psi (t_\e,x_\e)  + \Psi (t_\e,x_\e) - \Psi(t_\e - h, x_\e) + \psi_\e (t_\e - h, x_\e)  - \psi_\e(t_\e -\e a, x_\e)
		\end{multline*}
		Thanks to the local maximum property, the sum of the four first terms is non-positive. Il follows that:
		\begin{align*}
		III  \leq & \int_{\frac h \e}^{1 + \frac t \e} \int_\RR \Phi \left(a \e \right) \omega(z) e^{|z|L}
		\exp \left( \frac{1}{\e} [ \Psi (t_\e,x_\e) - \Psi(t_\e - h, x_\e) ] \right) \\
		& \quad \cdot \exp \left( \frac 1 \e \left[\psi_\e(t_\e - h ,x_\e) - \psi_\e(t_\e -\e a, x_\e )\right] \right) \,{\rm d}z\,{\rm d}a.
		\end{align*}
		Since $\partial_t \Psi \leq -\nu$ over $B( (t_\e, x_\e), h)$, we can bound $III$ as follows:
		\begin{align*}
		III  & \leq \int_{\frac h \e}^{1 + \frac t \e} \int_\RR \Phi \left(a \right) \omega(z) e^{|z|L} e^{-\frac {\nu h} \e} 
		\exp \left( \frac 1 \e \left[\psi_\e(t_\e - h ,x_\e) - \psi_\e(t_\e -\e a, x_\e )\right] \right) \,{\rm d}z\,{\rm d}a \\
		& \leq \int_{0}^{1 + \frac {t-h} \e} \int_\RR \Phi \left(a + \frac h \e \right) \omega(z) e^{|z|L} e^{-\frac {\nu h} \e} 
		\exp \left( \frac 1 \e \left[\psi_\e(t_\e - h ,x_\e) - \psi_\e(t_\e - h - \e a, x_\e )\right] \right) \,{\rm d}z\,{\rm d}a.
		\end{align*}
		Since equation~\eqref{CGM_eq_phi} is autonomous, we derive the following estimate on the time span $[t-h-l, t-h]$ in the same way we derived that of Lemma~~\ref{CGM_lem_m_e} on $[0,T]$:
		%
		\begin{equation}
		\exp \left( \frac 1 \e \left[\psi_\e(t_\e - h ,x_\e) - \psi_\e(t_\e - h - l, x_\e )\right] \right) \leq 1 + \frac{l}{\e} + K \e \left[1 + \frac{l}{\e}\right]^{1+\mu}.
		\end{equation}
		It follows that 
		\[
		III \leq  \int_0^{1 + \frac{t - h }{\e}} \exp \left( -\frac{\nu h}{\e} \right)  \Phi \left(a + \frac h \e \right) \left[ 1 + a + 
		K \e \left(1 + a\right)^{1+\mu} \right] \,{\rm d}a \int_\RR \omega(z) e^{|z|L} \,{\rm d}z.
		\]
		Since $\omega$ is a Gaussian distribution, the integral in $z$ (right factor) is finite. Since $\nu h > 0$, $t$ is bounded and $\Phi$ is algebraic, the integral in $a$ (left factor) converges to $0$ as $\e \to 0$.

		Passing to the limit $\e \to 0$ in~\eqref{CGM_eq_I_II} now gives us:
		\[
		1 - \delta \leq \int_0^\infty \int_\RR \Phi (a) \omega(z) \exp \left( [a \partial_t \Psi + z \partial_x \Psi](t_0, x_0)\right) \,{\rm d}z \,{\rm d}a.
		\]
		By taking the limit when $\delta \to 0$ we recover:
		\[
		1 \leq \int_0^\infty \Phi (a) \exp \left( a \partial_t \Psi(t_0,x_0) \right)  \,{\rm d}a  \int_\RR  \omega(z)  \exp \left( z \partial_x \Psi (t_0, x_0) \right) \,{\rm d}z.
		\]
		Therefore, $\psi_0$ is a viscosity supersolution of~\eqref{CGM_eq_limiting_HJ}. 
	\end{proof}

	\begin{proof}[Proof of Theorem~\ref{CGM_thm_viscosity_solution}]~
		Propositions \ref{CGM_prop_viscosity_subsolution} and \ref{CGM_prop_viscosity_supersolution} prove $\psi_0$ is a viscosity solution of the Hamilton-Jacobi equation~\eqref{CGM_eq_limiting_HJ}. Since $\psi_0$ is bounded below and Lipschitz continuous, and the Hamiltonian $H$ satisfies the pertinent hypotheses, Theorem~\ref{CGM_thm_uniqueness} proves that $\psi_0$ is the unique viscosity solution of~\eqref{CGM_eq_limiting_HJ}. Local compactness of $(\psi_\e)_\e$ and standard Hausdorff separation arguments prove that the whole sequence $\psi_\e$ tends to $\psi_0$.
	\end{proof}

	\begin{cor} \label{cor_exp_tails}
		Assume Hypothesis~\ref{CGM_hyp_IC} and replace Hypothesis~\ref{CGM_hyp_omega} by the following.
		Let $\omega$ be an isotropic multivariate continuous probability distribution of mean $0$ such that, for some positive $\delta$,
		\[
			\int_{\RR^d} \omega (z) \exp \left( (C_L + \delta) |z| \right) \,{\rm d}z < \infty,
		\]
		where $C_L$ is the Lipschitz constant in space of the initial condition introduced in Hypothesis~\ref{CGM_hyp_IC}.\ref{CGM_hyp_IC_regularity}.
		Let $\beta(a) = \mu / (1+a)$ with $0<\mu<1$ as previously.
		
		Then $\psi_\e \xrightarrow[\e \to 0]{L^\infty_{\rm{loc}}} \psi_0$, which is the unique viscosity solution of the limiting Hamilton-Jacobi equation~\eqref{CGM_eq_limiting_HJ} with initial condition $v(x)$ among the class of bounded below, $C_L$-Lipschitz continuous functions.
	\end{cor}
\begin{proof}
	Let us sum up step by step the sufficient changes to our proofs.
	\begin{itemize}
		\item Proposition~\ref{CGM_prop_H} is modified as follows:
		\begin{itemize}		
			\item The Hamiltonian function $H$ is well defined in~\eqref{CGM_eq_def_H}, albeit only over an open set containing strictly $B_{\RR^d}(0,C_L)$.
			\item The Hamiltonian satisfies an inequality similar to $H(p) \gtrsim 1+|p|^2$ over $B_{\RR^d}(0,C_L)$ and can be modified over $\RR^d \setminus B_{\RR^d}(0,C_L)$ into $\tilde H$, with $\restriction{\tilde H}{B_{\RR^d}(0,C_L)} = \restriction{H}{B_{\RR^d}(0,C_L)}$, so as to preserve coercivity. We will later prove that the family $(\psi_\e)_\e$ is uniformly Lipschitz equicontinuous in space with constant $C_L$, so the modification of $H$ does not affect the modified Hamilton-Jacobi equation.
			\item The proof of convexity holds.
		\end{itemize}
	\item Hence the hypotheses of Theorem~\ref{CGM_thm_uniqueness} hold: the modified limiting Hamilton-Jacobi equation
	\[
	\partial_t \psi_0(t,x) + \tilde H (\nabla_x \psi_0) (t,x) = 0
	\] admits a unique solution among the class of bounded below, Lipschitz continuous functions.\\
	\item The bounds of Theorem~\ref{CGM_thm_bounds} suffer the following alterations:
	\begin{itemize}
		\item[\eqref{CGM_eq_LBpsi_simple}] The proof of the $L^\infty$ lower bound is unaffected.
		\item[\eqref{CGM_eq_UBpsi_simple}] The constant $C$ appearing in the equation is modified, but it remains a finite, positive constant.
		\item[\eqref{CGM_eq_Bpsi_x_simple}] The proof of Lipschitz continuity in space is unaffected.
		\item[\eqref{CGM_eq_Bpsi_t_simple}] The Lipschitz bound in time maintains the same expression and is still finite.
	\end{itemize}

	\item Viscosity limit procedure, Section~\ref{CGM_sec_viscosity_solution}:
	\begin{itemize}
		\item Proposition~\ref{CGM_prop_viscosity_subsolution} (viscosity subsolution) is unaffected.
		\item As for the viscosity supersolution procedure, Lemma~\ref{CGM_lem_B_e} is unaffected, and the analytical expressions in Lemmas~\ref{CGM_lem_A_e} and~\ref{CGM_lem_m_e}, as well as in Proposition~\ref{CGM_prop_B_e_to_0} remain the same (with different constants, as defined). The symmetry of $\omega$ plays a role in Proposition~\ref{CGM_prop_viscosity_supersolution}, as does the boundedness of $\int_\RR \omega(z) \exp(C_L |z|) \,{\rm d}z$
		in the proofs of the bounds on the expressions $II$ (small ages, large spaces) and $III$ (large ages), since the Lipschitz continuity constant in space of $\psi_\e$ is $C_L$. 
	\end{itemize}
	\end{itemize}
It follows that $\psi_\e$ converges to $\psi_0$, the unique bounded below, Lipschitz continuous viscosity solution of $\partial_t u + \tilde H (\nabla_x u) = 0$. By construction, $\psi_0$ is a viscosity solution of $\partial_t u + H (\nabla_x u) = 0$, and it is the unique solution when restricting to the class of bounded below, $C_L$-Lipschitz continuous functions.
\end{proof}


	\section{Appendix}\label{CGM_sec_perspectives}
	
	There are three main aspects we would like to discuss in this \cred{appendix}. 	
	First we briefly discuss  the technical motivation for our choices of hypotheses and proof strategies, and give a synthetic presentation of the main ideas behind our work.
	Second, we will support and elaborate on the claim we make, that equation~\eqref{CGM_eq_limiting_HJ} is the same as the limiting Hamilton-Jacobi equation derived after renormalising $n$ by a non-stationary measure inspired by~\cite{BLM} that approaches a meaningful self-similar profile. Third, we will discuss a setting in which the jump rate $\beta$ depends not only on age but also on space. For the sake of simplicity, the two last parts are presented in dimension $d = 1$.

%

\subsection{Motivation, main ideas, and difficulties}\label{ssec_motivations}

Usually, similar limit problems for which the limit equation is averaged with respect to the fast variable (age $a$) are handled with the perturbed test function method introduced in~\cite{Evans_perturbed_test}, see for instance \cite{ Bouin2012, Bouin-Mirrahimi, Caillerie-CRAS}. However, in our setting the perturbed function would be naturally unbounded. Here, we bypass this issue by working directly on the boundary value of our solution~\eqref{CGM_eq_n}. Namely, we reduce the solution $\phi_\e(t,x,a)$ to the knowledge of $\psi_\e(t,x) = \phi_\e(t,x,0)$. Note that the reconstruction of $\phi_\e$ from $\psi_\e$ along characteristic lines makes the problem non local in time. That is the first main idea in this work. The two following subsections describe the two major difficulties that we have encountered.

\subsubsection{Corrected maximum principle} \label{CGM_ss_space_homo}

While defining the waiting time distribution $\Phi$ in~\eqref{CGM_eq_def_captial_Phi} in the model description subsection~\ref{CGM_ss_model}, we noted that 
the mean residence time of particles $\int_0^\infty a \Phi(a)\, \,{\rm d}a$ is infinite in the subdiffusive case $\beta (a) = \mu / (1+a)$ with $\mu \in (0,1)$. 

It is equivalent to say that the stationary distribution  $N_\infty(a) = \exp \left(- \int_{0}^a \beta(s) \,{\rm d}s \right)$ is not an integrable function. 
A classical result based on the maximum principle states that $L^\infty$ bounds on the ratio $n(t,\cdot) / N_\infty$ are propagated if true at time $t=0$, for the space homogeneous problem. This provides fruitful estimates if $N_\infty$ is integrable. Otherwise it is useless.  

Besides, it is well known that the long-time asymptotics of the space homogeneous problem follow a self-similar scaling \cite{Feller}: 
\begin{equation} \label{CGM_eq_self_sim_profile}
n (t,a) \simeq \frac{1}{1+t} W_\infty \left( \frac{a}{1+t} \right)\,,
\end{equation}
where $W_\infty $ is the Dynkin-Lamperti arc-sine law.
The precise description of the intermediate asymptotics (\emph{i.e.} the estimate of the distance between the solution at time $t$ and the self-similar profile) was the subject of~\cite{BLM}. As a side remark, this makes us expect $\ln n$ to be increasing in time as $\ln(1+t)$, leading to the apparition of a logarithmic term in time in any $L^\infty$ estimate of $\ln n$ over a compact time interval, as it is the case in equation~\eqref{CGM_eq_UB_psi_e_brut}. 


\subsubsection{Contribution of initial distribution at time $t$} \label{CGM_ssec_persistence}

The second main difficulty we have tackled appears in~Section~\ref{CGM_ssec_visc_supersol}, in the proof that the limit $\psi_0$ of a subsequence of $\psi_\e$ is a viscosity super-solution of the limiting Hamilton-Jacobi equation~\eqref{CGM_eq_limiting_HJ}. It stems from the long persistence of the initial condition in the renewal flux term.

In our proof, the renewal flux term at age $a=0$ is split into \cred{two} relative contributions: that of the particles which have already jumped before, and that of the particles which have never jumped before (respectively, $\mathfrak{A}_\e$ and $\mathfrak{B}_\e$ in equation~\eqref{CGM_eq_def_A_e}). We expect the latter to disappear in the limit $\e\to 0$. However, due to the heavy tail of the waiting time distribution, it may have a relatively high contribution, that we must bound from above.
We solve this issue by a refined estimate of the relative contribution which expresses an anomalous exponent:
\begin{equation}
\mathfrak{B}_\e(t,x)   \lesssim \left (\dfrac{\e}{t}\right )^{1+\mu} \dfrac{T}{\e}.
\end{equation}

\subsubsection{Comments on the choice of initial conditions}\label{ssec_IC_comments}

\paragraph{Lack of compatibility of the initial condition}
	The initial condition that we take is smooth enough in $x$ (Lipschitz continuous). 
	However, we do not require for it to be compatible in the sense that the influx relation at age $a=0$ is satisfied at time $t=0$ in~\eqref{CGM_eq_n}.
	As a consequence, we allow discontinuities along $t = \varepsilon a$. This means that in general, we may have:
	\begin{equation}
	n_\e (0, x, 0) = \int_0^1 \int_{\RR^d} \omega(z) \Phi(a) n_\e^0(x-\e z, a) \,{\rm d}z\,{\rm d}a \neq n_\e^0(x,0).
	\end{equation}
	Such compatibility is assumed in \cite[Chapter~3.4]{Perthame_TEB} to infer regularity with respect to $a$ in the space-homogeneous setting. Such regularity is not required in the present contribution.

	\paragraph{Semi-concavity regularity} \label{CGM_rem_SC}
	The second strong assumption is the semi-concavity of the initial condition. It is required  in order to handle the relative contribution to the boundary renewal term of what remains from the initial data at time $t$ (particles that are jumping for the first time). This implies local Lipschitz continuity as a by-product. Therefore, we have deemed reasonable to assume global Lipschitz continuity. The latter assumption is also in accordance with the uniqueness result that we use (Theorem \ref{CGM_thm_uniqueness}). 
	This hinders the use of half-relaxed limits~\cite{Barles, 2011_ABIL}, which have been designed to bypass Lipschitz estimates.
	\subsection{Renormalising by a non-stationary measure}\label{CGM_ssec_instationary}
	
	The idea of renormalising the solution of a kinetic equation by a stationary measure and studying some multiplicative perturbation term is classical. However, as has been shown in~\cite{BLM}, it cannot be applied here in a straightforward way 
	because, were a steady state to exist in self-similar variables for our equation, it would be infinite at age $0$, rendering the boundary condition a meaningless ``$\infty = \infty$'' equality. 
	Let us attempt to remain as close as possible to the underlying principle by using a function that corresponds to the pseudo-equilibrium of~\cite{BLM}.
	
	For any $t > 0$ and $0 < a < 1+t$, let
	\begin{equation}
	N(t,a) = (1+a)^{-\mu}(1+t-a)^{\mu-1}.
	\end{equation}
	We also set, for any $x \in \RR$, $t > 0$ and $0 < a < 1+t$:
	\begin{equation}
	u(t,x,a) = \frac {n(t,x,a)}{N(t,a)}.
	\end{equation}
	We define the following measure, for $t > 0$ and $0 < a < 1+t$:
	\begin{equation}
	\nu_t (a) = \beta(a) \frac {N(t,a)}{N(t,0)} = \frac{\mu (1+t)^{1-\mu}}{(1+a)^{1+\mu}(1+t-a)^{1-\mu}}.
	\end{equation}
	Direct computation gives us
	\[
	\partial_t \ln N + \partial_a \ln N + \beta(a) = 0,
	\]
	which is also satisfied by $n$. Hence, $u$ satisfies:
	\begin{equation}
	\label{CGM_eq_u}
	\left\{ \begin{array}{l} 
	\partial_t u(t,x,a) + \partial_a u(t,x,a) = 0\, , \qquad \qquad \qquad t \geq0, \quad a > 0\, , \quad x\in \mathbb{R}\smallskip\\
	u(t,x,0) = \int_0^{1+t} \int_\RR \nu_t(a) \omega(x-x') u(t,x',a) \,{\rm d}x' \,{\rm d}a' \smallskip\\
	u(0,x,a) = u^0(x,a) = n^0(x,a) N(0,a) \qquad \qquad {\text{with }} \quad {\rm supp} (u^0(x,\cdot)) = [0,1].
	\end{array}\right.
	\end{equation}

	Let us take a hyperbolic time - space scaling and a Hopf-Cole transform:
	\begin{equation}
	u_\e (t,x,a) = u \left( \frac t \e, \frac x \e, a \right) = \exp\left(-\frac 1 \e \tilde\varphi_\e^0 (t,x,a)\right).
	\end{equation}
	Characteristic flow of \eqref{CGM_eq_u} leads us to define:   
	\begin{equation}
	\tilde\varphi_\e (t,x,a) =
	\begin{cases} 
	\tilde\varphi_\e^0 (x, a - t/ \e), & \qquad a > t/\e \smallskip\\
	\tilde\psi_\e (t - \e a, x), & \qquad a \leq t / \e.
	\end{cases}
	\end{equation}
	Let us also set, in agreement with the Ansatz~\eqref{CGM_eq_Ansatz} in Hypothesis~\ref{CGM_hyp_IC}:
	\begin{equation}
	\begin{aligned}
	\tilde\varphi_\e^0 (x,a) & = v(x) + \e \xi (x,a) + \chi_{[0,1]}(a) \smallskip\\
	& = v(x) + \e \left[ \eta (x,a) - (1+\mu) \ln(1+a) - (1-\mu) \ln (1-a) \right] + \chi_{[0,1]}(a),
	\end{aligned}
	\end{equation}
	where $\chi_A$ is worth $0$ over the set $A$ and $+\infty$ outside of $A$.
	
	With the previous definitions, $\tilde\varphi_\e$ satisfies the following equation, which is analogous to~\eqref{CGM_eq_eq}:
	\begin{equation} \label{CGM_eq_eq_N} 
	\begin{aligned}
	1 = & \int_0^{t / \e} \int_\RR \exp\left(\frac 1 \e \left[ \tilde\psi_\e (t,x) - \tilde\psi_\e(t-\e a, x-\e z) \right]\right) \nu_{t/\e}(a) \omega(z) \,{\rm d}z\,{\rm d}a \\
	& + \int_0^1\int_\RR \exp\left(\frac 1 \e \left[ \tilde\psi_\e (t,x) - \tilde\varphi_\e^0 (x-\e z, a) \right]\right) \nu_{t/\e} (a+t/\e) \omega(z) \,{\rm d}z\,{\rm d}a.
	\end{aligned}
	\end{equation}
	
	\begin{rem}
		For any positive $t$,
		\[
		\int_0^{1+ t / \e} \nu_{t/\e}(a)\,{\rm d}a= \frac{1}{1 + \frac{1}{1+t/\e}} \xrightarrow[\e \to 0]{} 1.
		\]
	\end{rem}
	Assuming sufficient regularity, \eqref{CGM_eq_eq_N} gives us:
	\[
	1 = \int_0^{t / \e} \int_\RR \exp\left(a \partial_t \tilde\psi_\e(t,x)\right) \exp\left(z \partial_x \tilde\psi_\e (t,x)\right) \exp\left(o(1)\right) \omega(z) \Phi(a) \left(\frac{1 + t/\e}{1 - a + t /\e}\right)^{1-\mu} \,{\rm d}z\,{\rm d}a.
	\]
	Hence the formal limit of~\eqref{CGM_eq_eq_N} is the same Hamilton-Jacobi equation as~\eqref{CGM_eq_limiting_HJ}:
	\[
	1 = \int_0^\infty \Phi(a) \exp\left(a \partial_t \tilde\psi_0 (t,x)\right) \,{\rm d}a \int_\RR \omega(z) \exp\left(z \partial_x \tilde\psi_0 (t,x)\right) \,{\rm d} z,
	\]
	with the same initial condition $v$.
	
	\begin{rem}
		In order to prove convergence of this newly defined $\tilde\psi_\e$ to $\tilde\psi_0$, solution of the limiting Hamilton-Jacobi equation, the computations required are more or less the same as those presented in this article, but they contain an additional term that must be estimated:
		\begin{align*}
		\frac 1 {2 + t / \e} = & \int_0^{t / \e} \int_\RR \left[ \exp\left(\frac 1 \e \left[ \tilde\psi_\e (t,x) - \tilde\psi_\e(t-\e a, x-\e z) \right]\right) - 1 \right] \nu_{t/\e}(a) \omega(z) \,{\rm d}z\,{\rm d}a \smallskip\\
		& + \int_0^1 \int_\RR \left[  \exp\left(\frac 1 \e \left[ \tilde\psi_\e (t,x) - \tilde\varphi_\e^0 (x-\e z, a) \right]\right)  - 1 \right] \nu_{t / \e}(a+t / \e) \omega(z) \,{\rm d}z\,{\rm d}a.
		\end{align*}
		This is due to the fact that $\nu_{t/\e}$ is not a probability measure over $[0, 1+t/\e]$. Since $\nu_{t/\e}$ does approach a probability measure for any $t>0$ as $\e \to 0$, this is not a major problem. 
	\end{rem}

	\subsection{Space-dependent jump rate}
	
	Our study, as briefly mentioned in the Introduction, has a biological motivation. The random motion we model takes place in cellular media in which heterogeneities are often prevalent. Hence the relevance of considering a space-dependant jump rate $\beta (x,a)$. There are different pertinent ways of defining the jump rate, depending on what we intend to model. Here, we will only consider the simple case of a slow space variation of the jump rate, in the sense that follows. We define
	\begin{equation} 
	\beta (x,a) = \frac {\mu(x)}{1+a},
	\end{equation}
	where $0 < \mu < 1$ is Lipschitz continuous, and consider the following problem:
	\begin{equation}
	\label{CGM_eq_n_e_beta_x}
	\left\{ \begin{aligned}
	& \partial_t n_\e(t,x,a) + \frac 1 \e \partial_a n_\e(t,x,a) + \frac 1 \e \beta(x,a) n_\e(t,x,a) = 0\, ,\quad t \geq0, \quad a > 0\, , \quad x\in \mathbb{R}\smallskip\\
	& n_\e(t,x,0) = \int_0^{1+ t / \e} \int_\RR \beta(x - \e z, a) \omega(z) n_\e(t,x - \e z ,a) \,{\rm d}z \,{\rm d}a \smallskip\\
	& n_\e(0,x,a) = n_\e^0(x,a) = n^0(x/\e, a).
	\end{aligned} \right.
	\end{equation}

	\begin{rem}
		It follows that $n(t,x,a) = n_\e (\e t, \e x,a)$ satisfies the problem below, with a jump rate that varies slowly in space:
		\begin{equation} \label{CGM_eq_P_e}
		\left\{ \begin{aligned}
		& \partial_t n (t,x,a) + \partial_a n(t,x,a) + \beta (\e x,a) n (t,x,a) = 0\, ,\quad t \geq0, \quad a > 0\, , \quad x\in \mathbb{R}\smallskip\\
		& n(t,x,0) = \int_0^{1+ t } \int_\RR \beta( \e x' , a) \omega(x - x') n (t, x', a) \,{\rm d}x' \,{\rm d}a \smallskip\\
		& n (0,x,a) = n^0(x,a).
		\end{aligned} \right.
		\end{equation}
	\end{rem}
	
	Since $\mu$ is Lipschitz continuous, 
	\begin{equation} \label{CGM_eq_beta_x_slow}
	\beta (x-\e z,a) = \frac {\mu(x) + \mathcal{O}(\e z)}{1+a} .
	\end{equation}
	The formulation of~\eqref{CGM_eq_n_e_beta_x} along characteristic lines allows us to recover, for $\psi_\e$ and $\phi_\e^0$ defined as in~\eqref{CGM_eq_carac},
	\begin{multline} \label{CGM_eq_eq_beta_x}
	1 =  \int_0^{t/\e} \int_\RR \omega(z) \Phi(x-\e z,a) \exp\left(\frac 1 \e \left[ \psi_\e (t,x) - \psi_\e(t-\e a, x - \e z) \right]\right)\,  \,{\rm d} z\,  \,{\rm d} a \\
	+ \int_{t/\e}^{1+t/\e} \int_\RR \omega(z) \Phi(x - \e z, a) \exp\left(\frac 1 \e \left[ \psi_\e (t,x) - \phi_\e^0(x - \e z, a - t / \e) \right] + \int_0^{a - t/\e} \beta(x, s)\,  \,{\rm d} s\right)\,  \,{\rm d} z\,  \,{\rm d} a,
	\end{multline}
	where
	\begin{equation}
	\Phi(x,a) = \beta(x,a) \exp \left( - \int_0^a \beta(x,s) {\rm d}s \right).
	\end{equation}
	Thanks to \eqref{CGM_eq_beta_x_slow} and since $\omega$ is a Gaussian, it follows that \eqref{CGM_eq_eq_beta_x} admits a formal limiting Hamilton-Jacobi equation, similar to the space-independent case \eqref{CGM_eq_limiting_HJ}. Here however, the Hamiltonian depends on space:
	\begin{equation}
	\label{CGM_eq_limiting_HJ_beta_x}
	1 = \int_0^\infty \Phi(x,a) \exp\left(a \partial_t \psi_0 (t,x)\right) \,{\rm d}a \int_\RR \omega(z) \exp\left(z \partial_x \psi_0 (t,x)\right) \,{\rm d} z .
	\end{equation}
	Yet again, that is a Hamilton-Jacobi equation, since it is equivalent to:
	\begin{equation}
	\label{CGM_HJ_2}
	\partial_t \psi_0 (t,x) + H(x, \partial_x \psi_0)(t,x) = 0 ,
	\end{equation}
	with $H$ defined as follows, where $\left(\hat\Phi(x,\cdot)\right)^{-1}$ is the inverse function of the Laplace transform of $\Phi(x,\cdot)$: 
	\begin{equation} \label{CGM_eq_def_H_x}
	H(x,p) = - \left(\hat\Phi(x,\cdot)\right)^{-1}\left( \frac 1 { \int_\RR \omega(z) \exp(z p) \,{\rm d}z } \right) .
	\end{equation}
	Passing to the limit rigorously is left for further work.

	\section*{Acknowledgements}
	
	This work was initiated within the framework of the LABEX MILYON (ANR-10-LABX-0070) of Universit\'e de Lyon, within the programme `Investissements d'Avenir' (ANR-11-IDEX-0007) operated by the French National Research Agency (ANR).
	
	This project has received funding from the European Research Council (ERC) under the European
	Union’s Horizon 2020 research and innovation programme (grant agreement No. 639638).
	
	The second author was supported by the ANR project KIBORD. (ANR-13-BS01-0004).
	
	During the project, the third author was affiliated to the Unit\'e de Math\'ematiques Pures et Appliqu\'ees (UMPA) UMR CNRS $5669$, of the \'Ecole Normale Sup\'erieure de Lyon and to the Project-Team Beagle of the Inria Rh\^one-Alpes.
	
	We wish to thank Thomas Lepoutre and Hugues Berry for many fruitful discussions.
	
	We wish to thank the anonymous reviewer for his or her thorough and useful work.



\bibliographystyle{plain}           		
\bibliography{bibliography}        	

%

\end{document}